\documentclass[preprint]{elsarticle}
\usepackage{amssymb}
\usepackage{lineno,hyperref}
\usepackage{textcomp}
\usepackage{amsmath}
\usepackage{lipsum}
\usepackage{pdflscape}
\usepackage{longtable,array}
\usepackage{rotating}

 \usepackage{graphics}
 \usepackage{graphicx}

\usepackage{amssymb}
\usepackage{amsmath}
\usepackage{amssymb}
\usepackage{float}
\usepackage[justification=centering]{caption}
\usepackage{subfig}
\usepackage{lineno}
\usepackage[margin=1in]{geometry}
\usepackage{subfloat}
\newtheorem{theorem}{Theorem}
\newdefinition{rmk}{Remark}

\modulolinenumbers[5]

\journal{Archive}

\begin{document}

\begin{frontmatter}

\title{Similarity reductions, new traveling wave solutions, conservation laws of (2+1)- dimensional Boiti-Leon-Pempinelli system}
\author[mysecondaryaddress]{Subhankar Sil}
\author[mysecondaryaddress]{T. Raja Sekhar\corref{mycorrespondingauthor}}
\cortext[mycorrespondingauthor]{Corresponding author}
\ead{trajasekhar@maths.iitkgp.ac.in
}

\address[mysecondaryaddress]{Department of Mathematics, Indian Institute of Technology Kharagpur, Kharagpur-2, India}

\begin{abstract}
In this article we obtain exact solutions of (2+1)-dimensional Boiti-Leon-Pempinelli system of nonlinear partial differential equations which describes the evolution of  horizontal velocity component of water waves propagating in two directions. We perform the Lie symmetry analysis to the given system and  construct one-dimensional optimal subalgebra which involves some arbitrary functions of spatial variables. Several new exact solutions are obtained by symmetry reduction using each of the optimal subalgebra. We then study the physical behavior of some  exact solutions by numerical simulations and observed many interesting phenomena such as traveling waves, lump type solitons, kink and anti-kink type solitons, breather solitons, singular kink type solitons and etc. We construct several conservation laws of the system by using multipliers method. As an application, we study the nonlocal conservation laws of the system by constructing potential systems and appending gauge constraints.
\end{abstract}


\end{frontmatter}
\section{Introduction}
A wide range of nonlinear physical phenomena in the vast areas of scientific disciplines are depicted by the nonlinear coupled partial differential equations (PDEs). Many significant phenomena in physics and engineering are represented by such nonlinear PDEs. These systems describe multiple behaviors in various fields such as mathematical physics, fluid dynamics, chemistry, condensed matter, biophysics, plasma physics, optical fibers, biology and other areas of engineering. The exact solutions of such system of nonlinear PDEs play an important role in nonlinear science, especially in nonlinear physics, since they can yield very much physical information and more insight into the physical aspects of the problem and thus lead to applications like understanding the behavior of the physics associated with the problem and also to test and analyze numerical schemes. The exact solutions of nonlinear PDEs are very interesting and popular area of research in nonlinear mathematical physics. However, no effective method has been proposed till date to derive the general solution of nonlinear PDEs; only special solutions  can be obtained by a few methods such as inverse scattering transformation, B$\ddot{a}$cklund transformation, Darboux transformation, Hirota's direct method, Painleve analysis, symmetry reductions,  variable separation approach,  homogeneous balance method, F-expansion method and etc.

Symmetry analysis \cite{bluman2010applications} is one of the most efficient tool and easy to implement when searching for some particular exact solutions to differential equations. A symmetry of system of PDEs is  one-parameter Lie group of transformations which leaves the given system invariant, or more preciously, a symmetry of PDE system leaves the solution manifold of that system invariant and it maps one solution to another solution of the given PDE system. Once one has determined the symmetry group of a system of differential equations, based on it variety of applications are available. One of the most important and useful application of symmetry method is to obtain systematically some classes of exact solutions \cite{olver1987group,sahoo2020optimal,satapathy2018optimal,yacsar2011symmetries,saha2020lie2}. A particular solution obtained from  symmetry group $G$ is group invariant solution corresponding to the group $G$.  For a given PDE system there may exist infinitely many particular solutions, so one needs to minimize the search for exact solutions. In this context, the concept of classification of optimal subalgebras was introduced by Ovsiannikov \cite{ovsiannikov2014group} where one needs to find a set of inequivalent subalgebras. Later, Olver \cite{olver2000applications} improved this method by introducing adjoint representations. Recently, many mathematicians \cite{sekhar2016group,cherniha2021complete,sil2020nonclassical,vaneeva2020generalization,benoudina2021lie,opanasenko2017group} contributed in this direction and obtained exact solutions of various physically relevant systems.

Conservation laws \cite{zhang2021symmetry,liu2020existence} describe many essential physical properties of a given PDE system and have also applications in existence, uniqueness and stability analysis for the development of numerical methods. Moreover, one can construct nonlocally related PDE systems of the original PDE system by introducing some potential (nonlocal) variables through conservation laws. A PDE may have more than one conservation law which arises by multiplying appropriate multipliers \cite{anco1997direct,anco2002direct,anco2002direct1} to the given PDE. Recently, Sil et al. \cite{sil2020nonlocal,sil2020nonlocally} applied direct multipliers method to construct conservation laws and applied them to study nonlocal symmetry analysis.

For the case where the associated Lie algebra is of infinite dimensional, the classification of optimal subalgebra is of special interest. Here one  obtains an infinite number of exact solutions as the corresponding symmetries involve some arbitrary functions of independent variables or dependent variables.

The (2+1)-dimensional Boiti-Leon-Pempinelli (BLP) system \cite{BLP}
\begin{eqnarray}\label{govv}
&&u_{ty}-(u^2-u_x)_{xy}-2v_{xxx}=0,\\
\nonumber&&v_t-v_{xx}-2uv_x=0
\end{eqnarray}
 is actually a  generalization of the (2+1)-dimensional sinh-Gordon equations. The Hamiltonian structure, Painlev$\acute{e}$ property \cite{mu2013localized}, Lax pair \cite{BLP}, and B$\ddot{a}$cklund transformation \cite{jiang2010solitons,zhao2017lie}  have been studied for BLP system \eqref{govv} and moreover various exact solutions \cite{lu2004soliton,ma2003diversity,huang2004exact}
were obtained by using tanh-coth method \cite{wazwaz2010}, CTE solvability \cite{wang2014}, improved projective equation approach and a linear variable separation approach \cite{yang2011}. Later Kumar and Kumar \cite{kumar2014} studied the BLP system \eqref{govv} in terms of Lie symmetry analysis and obtained one family of solutions consisting of arbitrary function.
Very recently Wang et. al \cite{wang20202+} proposed the modified BLP system
\begin{eqnarray}
\label{gov}R_1:&&u_{ty}=a(u^2-u_x)_{xy}+bv_{xxx}\\
\nonumber R_2:&&v_t=cv_{xx}+guv_x
\end{eqnarray}
by introducing some new parameters $a,b,c$ and $g$ where $a, b, c$ and $g$ are real numbers. It describes the evolution of the horizontal velocity component of water waves propagating in $x$ and $y$ directions in an infinite narrow channel of constant depth.  In fact, it demonstrates the evolution of the horizontal component of the velocity of water waves propagating through an infinite narrow channel that maintains constant depth of the $x-y$ plane. Here $t$ denotes time, $x,y$ represents spatial variables, the dependent variables $u$ and $v$ demonstrate the velocity component in $x$ and $y$ directions respectively.
In \cite{wang20202+}, the authors obtained only the stationary domain walls solution of \eqref{gov}. There is a major research gap in the direction of obtaining exact solutions of the BLP system \eqref{gov}. Since Lie symmetry analysis is the most powerful tool to construct exact solutions of nonlinear system of PDEs, therefore our aim is to study the BLP system \eqref{gov} by means of Lie symmetry analysis and obtain several new exact solutions by constructing set of optimal subalgebras. Moreover, we construct conservation laws of the BLP system \eqref{gov} by direct multiplier method. The outline of our work is as follows:\par
In section \ref{s2}, we apply the Lie symmetry method to the BLP system \eqref{gov} and compute the infinite dimensional Lie algebra. We perform the optimal classification of one-dimensional subalgebra consisting of arbitrary functions in section \ref{s3}. section \ref{s4} deals with obtaining several new exact solutions systematically from each subalgebra which are reported first time in the literature and also discuss the physical significance of the solution profiles geometrically. We obtain several new traveling wave solutions  these indicate various important physical properties in section \ref{s5}.  In section \ref{s6}, we construct conserved vectors of the BLP system \eqref{gov}. We study the nonlocal conservation laws of system \eqref{gov} as an application of those conserved vectors in section \ref{s7}. Finally we provide concluding remarks in section \ref{s8}.

\section{Lie symmetry analysis}\label{s2}
We apply the Lie symmetry analysis to the given system \eqref{gov}. Let us consider a one-parameter($\epsilon$) infinitesimal Lie group of point transformations of the form
\begin{eqnarray*}
&&t^\star=t+\epsilon \tau(t,x,y,u,v)+O(\epsilon^2),\\
&&x^\star=x+\epsilon \xi(t,x,y,u,v)+O(\epsilon^2),\\
&&y^\star=y+\epsilon \pi(t,x,y,u,v)+O(\epsilon^2),\\
&&u^\star=u+\epsilon \eta(t,x,y,u,v)+O(\epsilon^2),\\
&&v^\star=v+\epsilon \phi(t,x,y,u,v)+O(\epsilon^2),
\end{eqnarray*}
where $\epsilon$ is a parameter and $\tau,\xi,\pi,\eta$ and $\phi$ are unknown infinitesimals which are to be determined. The associated Lie symmetry generator takes the form
$$\Psi=\tau(t,x,y,u,v)\frac{\partial}{\partial t}+\xi(t,x,y,u,v)\frac{\partial}{\partial x}+\pi(t,x,y,u,v)\frac{\partial}{\partial y}+\eta(t,x,y,u,v)\frac{\partial}{\partial u}+\phi(t,x,y,u,v)\frac{\partial}{\partial v}.$$
Suppose $\Psi^{(3)}$ is a 3rd prolongation of $\Psi$, then the symmetry determining equations are
\begin{eqnarray*}
&&\Psi^{(3)}(R_1)|_{\{R_1=0,R_2=0\}}=0,\\
&&\Psi^{(3)}(R_2)|_{\{R_1=0,R_2=0\}}=0,
\end{eqnarray*}
where the given PDE system is denoted as $R_1=0,R_2=0$ which results an overdetermined linear system of PDEs. After solving the determining system we obtain the unknown infinitesimals as
\begin{eqnarray*}
&&\tau=c_3-2c_4t,\\
&&\xi=c_1-2c_4x,\\
&&\pi=c_2-F_1(y),\\
&&\eta=c_4u,\\
&&\phi=F_2(y)+vF_1'(y)
\end{eqnarray*}
where $c_1,c_2,c_3$ and $c_4$ are arbitrary constants while $F_1(y)$ and $F_2(y)$ are arbitrary functions of $y$. Consequently the Lie symmetry generators are listed as
\begin{eqnarray*}
&&X_1=\frac{\partial}{\partial x},\\
&&X_2=\frac{\partial}{\partial y},\\
&&X_3=\frac{\partial}{\partial t},\\
&&X_4=-2t\frac{\partial}{\partial t}-x\frac{\partial}{\partial x}+u\frac{\partial}{\partial u},\\
&&X_5=F_2(y)\frac{\partial}{\partial v},\\
&&X_6=-F_1(y)\frac{\partial}{\partial y}+vF_1'(y)\frac{\partial}{\partial v}.
\end{eqnarray*}

Now, the corresponding one-parameter Lie group of transformations can be obtained by solving the following initial value problem:
\begin{eqnarray*}
&&\frac{d\tau}{d\epsilon}=\tau(t,x,y,u,v),\\
&&\frac{d\xi}{d\epsilon}=\xi(t,x,y,u,v),\\
&&\frac{d\pi}{d\epsilon}=\pi(t,x,y,u,v),\\
&&\frac{d\eta}{d\epsilon}=\eta(t,x,y,u,v),\\
&&\frac{d\phi}{d\epsilon}=\phi(t,x,y,u,v),
\end{eqnarray*}
with the initial data $\tau=t,\xi=x,\pi=y,\eta=u$ and $\phi=v$ when $\epsilon=0.$ On solving the above system of ODEs for each infinitesimal transformation, we obtain the corresponding one-parameter Lie group of transformations
\begin{eqnarray*}
&&G_1:~~(t^\star,x^\star,y^\star,u^\star,v^\star)=(t,x+\epsilon,y,u,v),\\
&&G_2:~~(t^\star,x^\star,y^\star,u^\star,v^\star)=(t,x,y+\epsilon,u,v),\\
&&G_3:~~(t^\star,x^\star,y^\star,u^\star,v^\star)=(t+\epsilon,x,y,u,v),\\
&&G_4:~~(t^\star,x^\star,y^\star,u^\star,v^\star)=\left(te^{-2\epsilon},xe^{-\epsilon},y,ue^{\epsilon},v\right).
\end{eqnarray*}

It is difficult to compute the group of transformations associated with the symmetries $X_5$ and $X_6$ since they involve some arbitrary functions and thus cannot be integrated. In the view of the above discussion on one-parameter Lie group of point transformations we have the following result:
\begin{theorem}
Let $u=g_1(t,x,y)$ and $v=g_2(t,x,y)$ be a solution of the given system \eqref{gov}. The group actions $G_1,G_2,G_3$ and $G_4$ acting on the solution surface $u-g_1=0,~v-g_2=0$ provide  one-parameter family of solutions
\begin{eqnarray*}
&&G_1: u=g_1(t,x-\epsilon,y),v=g_2(t,x-\epsilon,y),\\
&&G_2:u=g_1(t,x,y-\epsilon),v=g_2(t,x,y-\epsilon),\\
&&G_3:u=g_1(t-\epsilon,x,y),v=g_2(t-\epsilon,x,y),\\
&&G_4:u=e^\epsilon g_1\left(te^{2\epsilon},xe^\epsilon,y\right),v=g_2\left(te^{2\epsilon},xe^\epsilon,y\right).
\end{eqnarray*}
\end{theorem}

\section{Classification of optimal subalgebras}\label{s3}
In this section, we  discuss the structure of the infinite dimensional Lie algebra $\mathfrak{L}$. The commutator table, corresponding to the symmetries $X_i~\text{for~~}i=1,...,6,$ is presented in the Table \ref{table1} where the entry in the $ij$-th position of the table is defined as
$$[X_i,X_j]=X_iX_j-X_jX_i~~\text{for}~i,j=1,...,6.$$

\begin{table}[h!]
\centering
\begin{tabular}{|c||c|c|c|c|c|c|}
  \hline
  $[X_i,X_j]$ & $X_1$ & $X_2$ & $X_3$ & $X_4$ & $X_5(F_2)$ & $X_6(F_1)$ \\\hline\hline
  $X_1$ & 0 & 0 & 0 & $-X_1$ & 0 & 0 \\\hline
  $X_2$ & 0 & 0 & 0 & 0 & $X_5(F_2'(y))$ & $X_6(F_1'(y))$ \\\hline
  $X_3$ & 0 & 0 & 0 & $-2X_3$ & 0 & 0 \\\hline
  $X_4$ & $X_1$ & 0 & $2X_3$ & 0 & 0 & 0 \\\hline
  $X_5(F_2)$ & 0 & $-X_5(F_2'(y))$ & 0 & 0 & 0 & $X_5(F_1F_2'+F_2F_1')$ \\\hline
  $X_6(F_1)$ & 0 & $-X_6(F_1'(y))$ & 0 & 0 & $-X_5(F_1F_2'+F_2F_1')$ & 0 \\
  \hline
\end{tabular}
\caption{Commutator table corresponding to the Lie algebra $\mathfrak{L}$}
\label{table1}
\end{table}

For the construction of inequivalent set of optimal subalgebras, we need to find the adjoint representation of the symmetries. The adjoint action on $\mathfrak{L}$ is defined by the adjoint operator as
$$\text{Ad}_{\exp(\epsilon X_i)}X_j=e^{-\epsilon X_i}X_je^{\epsilon X_i},$$
where $\epsilon$ being the small parameter.

\begin{table}[h!]
\centering
\begin{tabular}{|c||c|c|c|c|c|c|}
  \hline
  $Ad_{e^{\epsilon X_i}}(X_j)$ & $X_1$ & $X_2$ & $X_3$ & $X_4$ & $X_5(F_2)$ & $X_6(F_1)$ \\\hline\hline
  $X_1$ &       $X_1$ &          $X_2$ &           $X_3$ &            $X_4+\epsilon X_1$ &        $X_5(F_2)$ &               $X_6(F_1)$ \\\hline
  $X_2$ &       $X_1$ &          $X_2$ &           $X_3$ &             $X_4$ &                 $X_5(F_2(y-\epsilon))$ &     $X_6(F_1(y-\epsilon))$\\\hline
  $X_3$ &       $X_1$ &          $X_2$ &           $X_3$ &           $X_4+2\epsilon X_3$  &      $X_5(F_2)$        &        $X_6(F_1)$ \\\hline
  $X_4$ & $e^{-\epsilon}X_1$  & $X_2$ &      $e^{-2\epsilon}X_3$ &    $X_4$ &               $X_5(F_2)$ &                   $X_6(F_1)$ \\\hline
  $X_5(F_2)$ &       $X_1$ &   $X_2+\epsilon X_5(F_2')$ &   $X_3$ &      $X_4$ &            $X_5(F_2)$ &          $X_6(F_1)-\epsilon X_5(F_1F_2'+F_2F_1')$ \\\hline
  $X_6(F_1)$ & $X_1$ &    $X_2+\epsilon X_6(F_1')$ &  $X_3$ &         $X_4$ &                    $X_5(F_2+\epsilon (F_1F_2)')$ &         $X_6(F_1)$ \\
  \hline
\end{tabular}
\caption{Adjoint table corresponding to the Lie algebra $\mathfrak{L}$}
\label{table2}
\end{table}
It can be defined as an infinite series form involving the Lie bracket which is given below:
$$\text{Ad}_{\exp(\epsilon X_i)}X_j=X_j-\epsilon[X_i,X_j]+\frac{\epsilon^2}{2}[X_i,[X_i,X_j]]-....$$
The adjoint actions are summarized in the Table \ref{table2}. It is to be noted that, due to  complexity of the calculations  we have used first two terms of the infinite series for computation of $\text{Ad}_{\exp(\epsilon X_6)}X_5$ without affecting the mathematical analysis for constructing the inequivalent optimal subalgebra. Now we  perform the classification of  inequivalent subalgebra using the adjoint table (refer Table \ref{table2}). First consider a general element
$$E=a_1X_1+a_2X_2+a_3X_3+a_4X_4+X_5(F_2)+X_6(F_1)$$
where $a_1,a_2,a_3$ and $a_4$ are constants. The basic idea \cite{Malek2018Lie} is to make $E$ more simpler element, say, $E'$ by choosing appropriate constants while using suitable adjoint actions. Suppose we apply the adjoint action of $X_i$ on $E$, then the updated element is of the form

\begin{eqnarray}\label{ad}
&&E'=\text{Ad}_{\exp(\epsilon X_i)}E=a_1'X_1+a_2'X_2+a_3'X_3+a_4'X_4+X_5+X_6,
\end{eqnarray}

where $a_i'$'$s$ are functions of $a_i$ and $\epsilon$.
 Here we perform the following cases:\par
$\textbf{Case-I:}$ We set $a_1\neq 0$ and let other constants be unrestricted. Without loss of generality we assume that $a_1=1.$ By choosing $X_i=X_1$ in \eqref{ad} we have
$$E'=(1+\epsilon a_4)X_1+a_2X_2+a_3X_3+a_4X_4+X_5(F_2)+X_6(F_1). $$
We cancel the $X_1$ term by choosing $\epsilon=-\frac{1}{a_4}$ and consequently we have
 $$E'=a_2X_2+a_3X_3+a_4X_4+X_5(F_2)+X_6(F_1).$$ Then we apply the adjoint action of $X_3$ on $E'$, which yields
 $E''=a_2X_2+(a_3+2\epsilon a_4)X_3+a_4X_4+X_5(F_2)+X_6(F_1).$ Again by choosing $\epsilon=-\frac{a_3}{2a_4}$, we are having with
 $$E''=a_2X_2+a_4X_4+X_5(F_2)+X_6(F_1).$$ Similarly applying the successive adjoint actions of $X_5,$  $X_6$ and canceling the $X_5,$ $X_6$ terms,  we have $E'''=a_2X_2+a_4X_4,$ where $a_2\in\{-1,0,1\}$ since generalized BLP system \eqref{gov} admits the discrete symmetries $(t,x,y,u,v)\rightarrow (t,x,-y,u,-v)$. Hence we assume $a_2=0,a_2=1$ which results the optimal set in this case as $E_1=X_2+a_4X_4$ and $E_2=X_4.$\par
$\textbf{Case-II:}$
\begin{table}[h]
\begin{tabular}{|c||c|c|c|c|c|c|}
  \hline
  $Ad_{e^{\epsilon X_i}}(E)$ & Coeff $X_1$ & Coeff $X_2$ & Coeff $X_3$ & Coeff $X_4$ & Coeff $X_5$ & Coeff $X_6$ \\\hline\hline
  $X_1$ & $a_1+\epsilon a_4$ & $a_2$ & $a_3$ & $a_4$ & $F_2$ & $F_1$ \\\hline
  $X_2$ & $a_1$ & $a_2$ & $a_3$ & $a_4$ & $F_2(y-\epsilon)$ & $F_1(y-\epsilon)$ \\\hline
  $X_3$ & $a_1$ & $a_2$ & $a_3+2\epsilon a_4$ & $a_4$ & $F_2$ & $F_1$ \\\hline
  $X_4$ & $a_1e^{-\epsilon}$ & $a_2$ & $a_3e^{-2\epsilon}$ & $a_4$ & $F_2$ & $F_1$ \\\hline
  $X_5$ & $a_1$ & $a_2$ & $a_3$ & $a_4$ & $F_2+\epsilon(a_2F_2'-F_1F_2'-F_2F_1')$ & $F_1$ \\\hline
  $X_6$ & $a_1$ & $a_2$ & $a_3$ & $a_4$ & $F_2+\epsilon (F_1F_2'+F_2F_1')$ & $F_1+a_2\epsilon F_1'$ \\
  \hline
\end{tabular}
\caption{Alternative way of representing adjoint table corresponding to the Lie algebra $\mathfrak{L}$}
\label{table3}
\end{table}

Now, let us assume the case $a_1=0,$ $a_2\neq 0$ and we consider the scaling of $E$ whenever needed. Further we assume that $a_2=1$ which leads to $E=X_2+a_3X_3+a_4X_4+X_5(F_2)+X_6(F_1).$ By choosing $X_i=X_5$ in \eqref{ad} we have
$$E'=X_2+a_3X_3+a_4X_4+X_5(F_2+\epsilon[F_2'-F_1F_2'-F_1'F_2])+X_6(F_1). $$ We eliminate the $X_5$ term by setting $F_2+\epsilon[F_2'-F_1F_2'-F_1'F_2]=0$ which yields
$E'=X_2+a_3X_3+a_4X_4+X_6(F_1)$. Now we apply the adjoint action of $X_6$ on $E'$ which provides
$$E''=X_2+a_3X_3+a_4X_4+X_6(F_1+\epsilon F_1').$$ Again we cancel the $X_6$ term by setting $F_1+\epsilon F_1'=0$ which in turn $E''=X_2+a_3X_3+a_4X_4.$ Once more applying the adjoint action of $X_3$ on $E''$ and canceling the $X_3$ term by proper choice of $\epsilon$, we are finally left with
$X_2+a_4X_4$ which is identical with $E_1$.\par


$\textbf{Case-III:}$
Now consider the case $a_1=0,a_2=0$ and $a_3=1.$ Here the general element $E=X_3+a_4X_4+X_5(F_2)+X_6(F_1).$ By choosing $X_i=X_3$ in \eqref{ad} we have
$E'=(1+2\epsilon a_4)X_3+a_4X_4+X_5(F_2)+X_6(F_1)$ which reduces to

$$E'=a_4X_4+X_5(F_2)+X_6(F_1)$$
after canceling the $X_3$ term by setting $\epsilon=-\frac{1}{2a_4}$. Again applying the adjoint action of $X_5$ on $E'$, it becomes $E''=a_4X_4+X_5(F_2+\epsilon[-F_1F_2'-F_1'F_2])+X_6(F_1)$. After canceling the $X_5$ term by choosing appropriate value  of $\epsilon$, finally we are left with another subalgebra $E_3=a_4X_4+X_6(F_1)$.  \par


$\textbf{Case-IV:}$ Here we consider $a_1=0,a_2=0,a_3=0$ and $a_4=1$. Take the general subalgebra element $E=X_4+X_5(F_2)+X_6(F_1)$ and applying the adjoint action of $X_5$ on $E$ and canceling out the $X_5$ term which leads to the subalgebra $X_4+X_6(F_1)$ which is equivalent to the subalgebra $E_3$.\par

$\textbf{Case-V:}$ We consider the case $a_1=0,a_2=0,a_3=0,a_4=0$ and $F_2\neq 0$. So here the general element is $E=X_5(F_2)+X_6(F_1)$. Applying the adjoint action on $E$ by $X_5$ and canceling the $X_5$ term with the possible course of action of $\epsilon$ we are left with $E_4=X_6(F_1)$.\par
From the above discussion we conclude the following result:
\begin{theorem}
The optimal system of one-dimensional subalgebras of the generalized BLP system \eqref{gov} consists of the following vector fields:
$$E_1=<X_2+a_4X_4>,~~E_2=<X_4>,~~E_3=<a_4X_4+X_6(F_1)>,~~E_4=<X_6(F_1)>.$$
\end{theorem}

\section{Similarity reductions and invariant solutions}\label{s4}
In this section, we obtain some group invariant solutions of the governing system \eqref{gov} by using each subalgebra in the optimal set.

\subsection{Reduction using $<X_2+a_4X_4>$}
In this case, the corresponding characteristic equations are
\begin{eqnarray*}
&&\frac{dt}{-2a_4t}=\frac{dx}{-a_4x}=\frac{dy}{1}=\frac{du}{a_4u}=\frac{dv}{0}.
\end{eqnarray*}
The similarity variables are $m=\frac{x}{\sqrt t},n=\frac{\ln(t)+2a_4y}{2a_4}$ and the corresponding similarity forms are
\begin{eqnarray*}
&&u(t,x,y)=\frac{U(m,n)}{\sqrt t},~~v(t,x,y)=V(m,n)
\end{eqnarray*}
where $U(m,n)$ and $V(m,n)$ are functions of $m$ and $n$ which are to be determined. Now, with this expression of $u$ and $v$, the governing system \eqref{gov} reduces to the following system with fewer independent variables
\begin{eqnarray}\label{red1}
&&a_4mU_{mn}-U_{nn}+a_4U_n+4aa_4U_mU_n-4aa_4UU_{mn}+2aa_4U_{mmn}-2ba_4V_{mmm}=0,\\
\nonumber&&a_4mV_m-V_n+2ca_4V_{mm}+2a_4gUV_m=0.
\end{eqnarray}
In general it is not feasible to solve so we apply the Lie symmetry technique to compute some exact solutions of \eqref{red1}. Using the Lie symmetry reduction method, we have the following ansatz:
$$U(m,n)=f_1(m),~~V(m,n)=\frac{pn}{s}+f_2(m)$$
where $p,s$ are arbitrary constants and $f_1(m),f_2(m)$ are unknown functions which are to be determined. With substitution of $U$ and $V$ in the above reduced system \eqref{red1} and solving the reduced system of ODEs we have
$$f_1(m)=\frac{1}{2}\frac{p-2sa_4cC_1-a_4m(C_1m+C_2)}{sa_4g(C_1m+C_2)},f_2(m)=\frac{1}{2}C_1m^2+C_2m+C_3$$
where $C_1,C_2$ and $C_3$ are integration constants. Thus, we have the  solution for the given system \eqref{gov} as follows:

\begin{figure}[h!]
  \centering
  \subfloat[]{\includegraphics[scale=0.28]{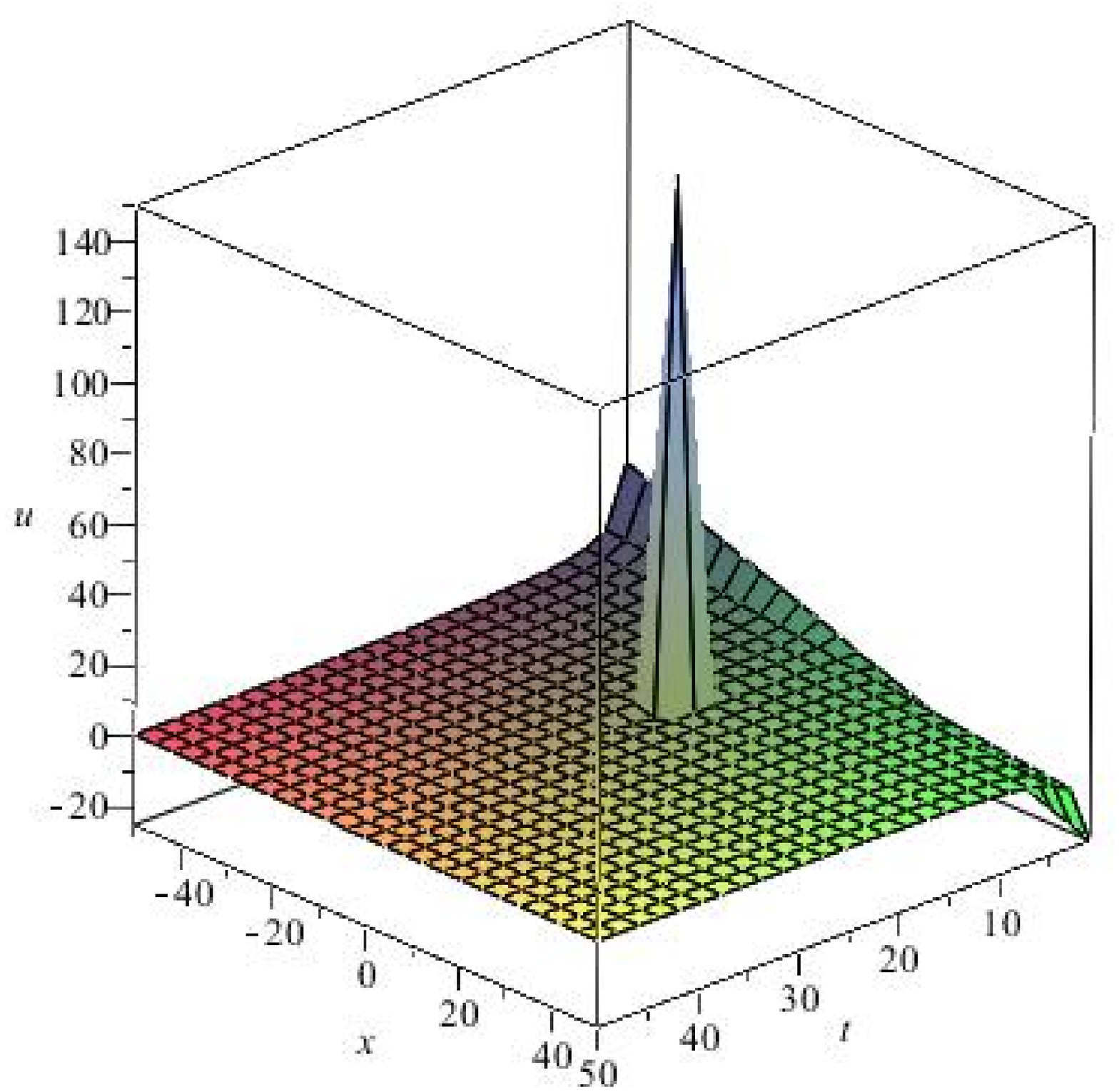}\label{sol1u2}}
  \\
  \subfloat[]{\includegraphics[scale=0.28]{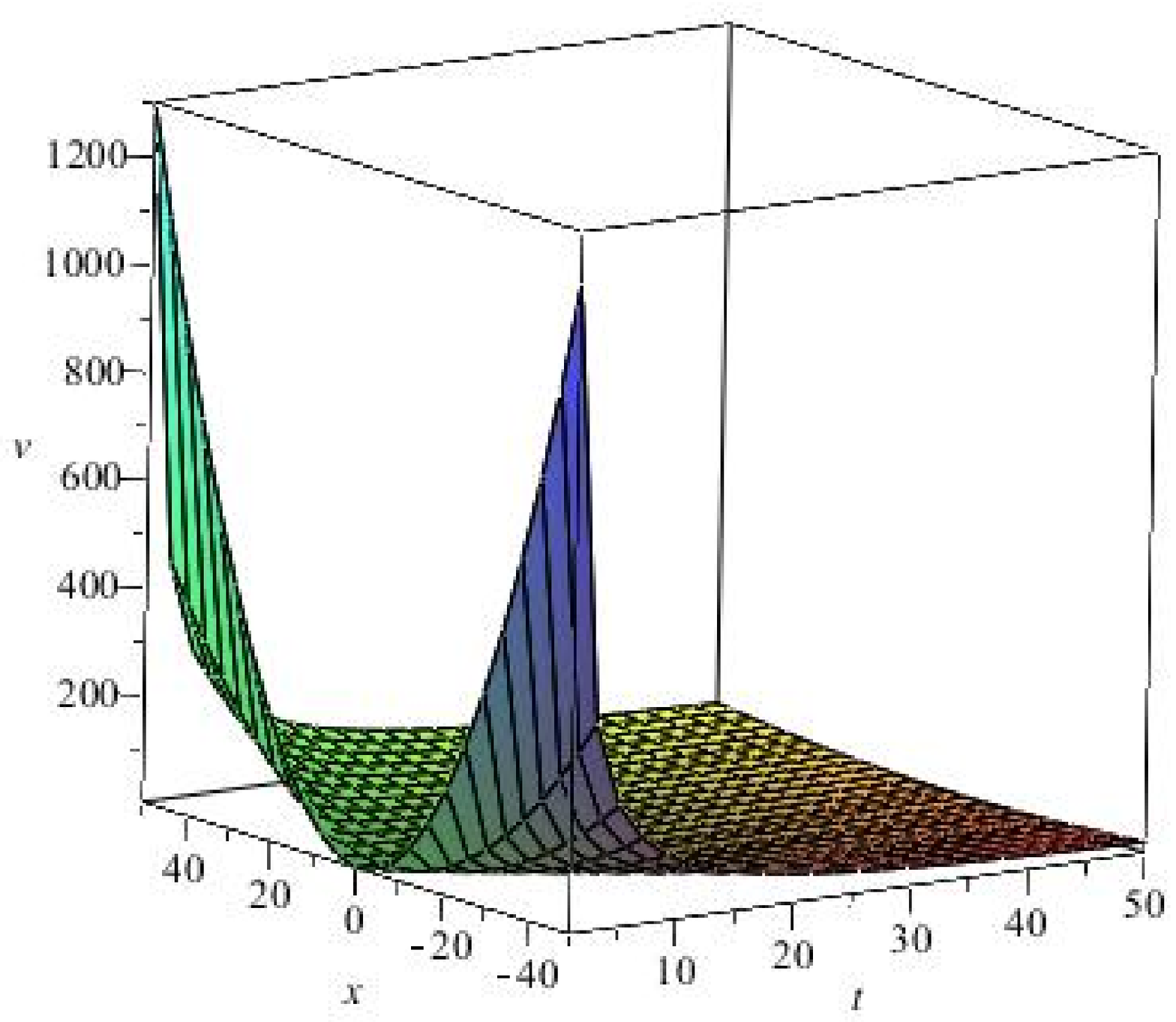}\label{sol1v2}}
  \subfloat[]{\includegraphics[scale=0.22]{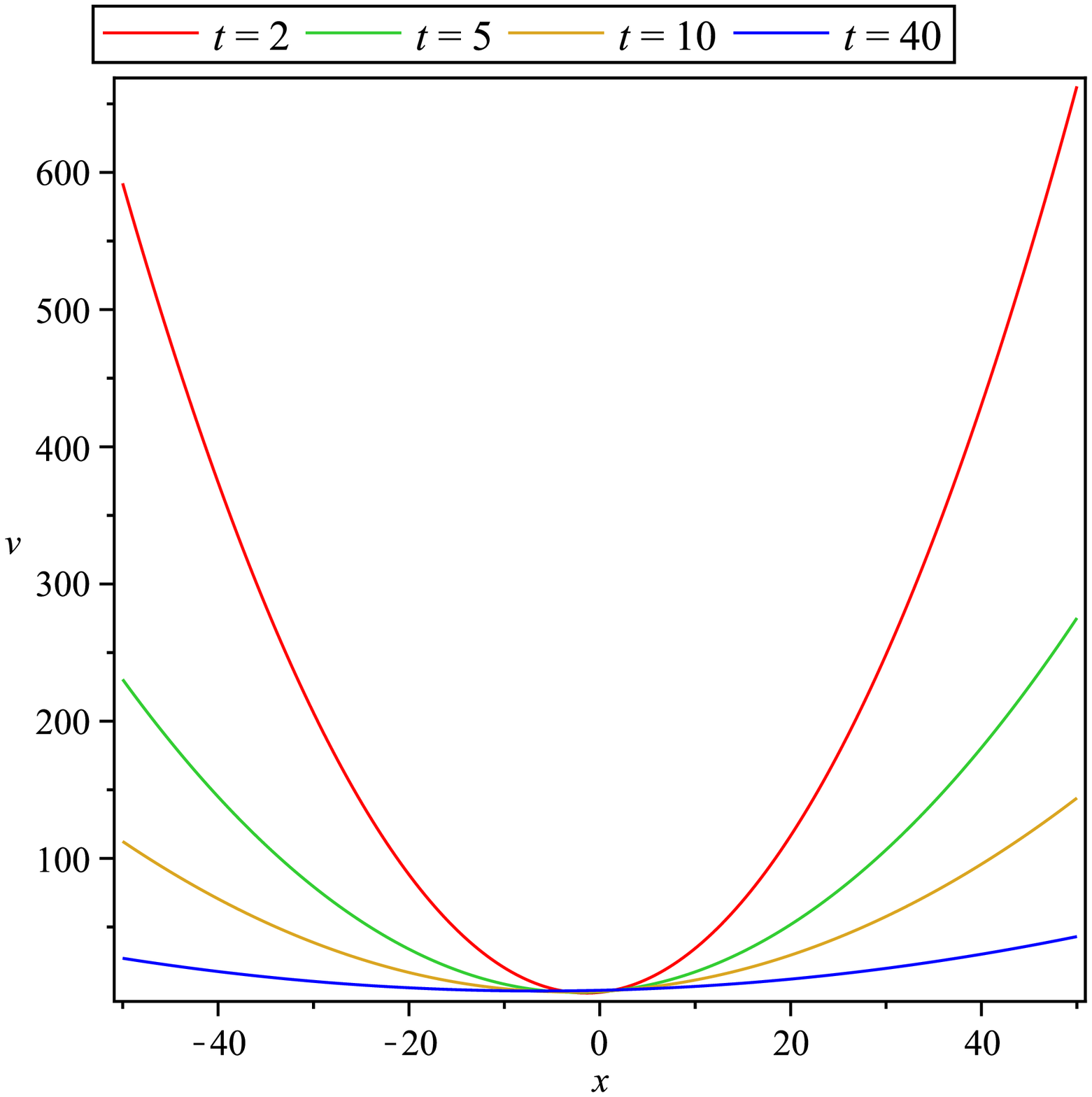}\label{sol1v22}}
  \caption{Solution profile of \eqref{gov} for a solution \eqref{sol1}: (a) 3d profile of $u_1$ (b) 3d profile of $v_1$ when $y=1$ (c) 2d profile of $v_1$ when $y=1$ for various values of $t$. }
  \label{fsol1}
  \end{figure}

\begin{eqnarray}\label{sol1}
&&u_1=-\frac{1}{2}\frac{a_4sC_1x^2+a_4sC_2x\sqrt t+(2cC_1a_4s-p)t}{a_4sgt(C_1x+C_2\sqrt t)},\\
\nonumber&&v_1=\frac{1}{2}\frac{a_4sC_1x^2+2a_4sC_2x\sqrt t+(2a_4C_3+p\ln(t)+2a_4py)t}{a_4s\sqrt t}.
\end{eqnarray}

The physical behavior of the solution profile \eqref{sol1} is illustrated in the Figure \ref{fsol1} by choosing $a_4=1,s=1,C_1=1,C_2=1,g=1$ and $p=1$. Figure \ref{sol1u2} represents single-lump soliton or 1-lump soliton for $u_1$. We illustrate the solution profile of $v_1$ in the Figure \ref{sol1v2} by fixing $y=1$. We noticed that (see, Figure \ref{sol1v2}) rapid increase of $v_1$ near the initial time when we approach far away from the origin that corresponds to two peaks. As time evolves, it gradually decreases uniformly when increasing the values of $x$.  We depict the 2-dimensional plot of $v_1$ by fixing $y=1$ with respect to $x$ for various values of $t$ in the Figure \ref{sol1v22}. It represents an upward parabola with vertex at the origin and as time evolves the parabola started to flatten and tends to a straight line.


By choosing the ansatz as $U(m,n)=f_3(n)$ and $V(m,n)=f_4(n)$ and plugging into the reduced system of PDEs \eqref{red1} we solve for the unknowns $f_3,f_4$ and obtain $$f_3(n)=C_4\exp{(a_4n)}+C_5,~f_4(n)=C_6,$$
where $C_4,C_5$ and $C_6$ are integration constants. So, we have the following exact solution for the given system \eqref{gov}
\begin{eqnarray}\label{sol2}
&&u_2=\frac{C_4\exp{\left[\frac{1}{2}(\ln(t)+2a_4y)\right]}+C_5}{\sqrt t},~~v_2=C_6.
\end{eqnarray}

\subsection{Reduction using $<X_4>$}
In this case, the similarity variables are $y$ and $m=\frac{x}{\sqrt t}$ which leads to the  invariant solution of the following form
$$u(t,x,y)=\frac{1}{\sqrt t}U\left(m,y\right),v(t,x,y)=V\left(m,y\right).$$

Substituting this form of $u$ and $v$ into the given system \eqref{gov}, we derive the reduced system
\begin{eqnarray}\label{red2}
&&mU_{my}+U_y+4aU_mU_y+4aUU_{my}-2aU_{mmy}+2bV_{mmm}=0,\\
\nonumber&&mV_m+2cV_{mm}+2gUV_m=0.
\end{eqnarray}
After solving the above system \eqref{red2}, we have the following solution
\begin{eqnarray*}
&&U(m,y)=\frac{1}{2g}\frac{C_7(-m^2-2c)-C_8m}{C_7m+C_8},\\
&&V(m,y)=\frac{C_7}{2}f_3(y)m^2+C_8f_3(y)m+g_3(y),
\end{eqnarray*}

\begin{figure}[h!]
  \centering
  \subfloat[]{\includegraphics[scale=0.28]{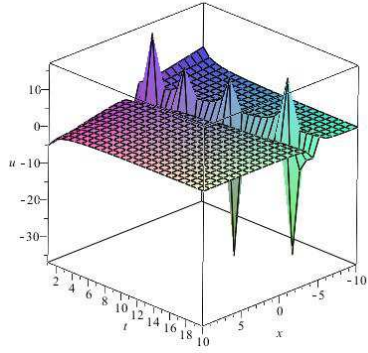}\label{sol3u2}}

  \subfloat[]{\includegraphics[scale=0.28]{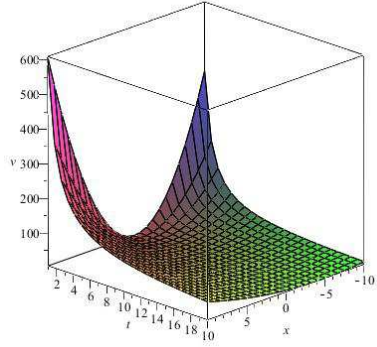}\label{sol3v2}}
  \subfloat[]{\includegraphics[scale=0.22]{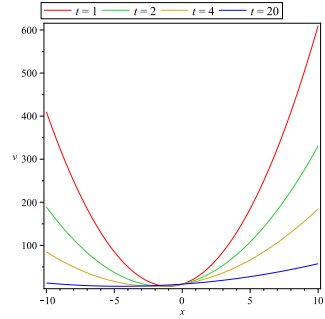}\label{sol3v22}}
  \caption{Solution profile of \eqref{gov} for a solution \eqref{sol3}: (a) 3d profile of $u_3$ (b) 3d profile of $v_3$ when $y=10$ (c) 2d profile of $v_3$ when $y=10$ for various values of $t$. }
  \label{fsol3}
  \end{figure}

where $C_7$ and $C_8$ are arbitrary constants and $f_3(y),g_3(y)$ are arbitrary functions of $y$, which in turn the solution of the given system \eqref{gov} as
\begin{eqnarray}\label{sol3}
&&u_3=-\frac{1}{2}\frac{C_7x^2\sqrt t+2cC_7t^{\frac{3}{2}}+C_8tx}{gt^{\frac{3}{2}}(C_7x+C_8\sqrt t)},\\
\nonumber&&v_3=\frac{1}{2}\frac{C_7x^2\sqrt tf_3(y)+2C_8xtf_3(y)+2t^{\frac{3}{2}}g_3(y)}{t^{\frac{3}{2}}}.
\end{eqnarray}
By choosing the parameters $c=1,g=1,C_7=1,C_8=1$ and considering $f_3(y)=y,g_3(y)=y$ we depict the solution profile which represents \eqref{sol3} (see, Figure \ref{fsol3}) for $u_3$ and $v_3$. The physical behavior of $u_3$ is demonstrated in the Figure \ref{sol3u2} which represents a multiple breather soliton. On the other hand, the 3-dimensional profile of $v_3$ is presented in the Figure \ref{sol3v2} by fixing $y=10$ and the 2-dimensional profile is illustrated in the Figure \ref{sol3v22} with respect to $x$ at various values of $t$ by fixing $y = 10$ which demonstrates an upward parabola.



We have another  solution of \eqref{red2} as
\begin{eqnarray*}
&&U(m,y)=-\frac{m}{4a}+\tanh(y-m),V(m,y)=f_5(y)
\end{eqnarray*}
where $f_5(y)$ is an arbitrary function of $y$, which yields a solution of the given system as
\begin{eqnarray}\label{sol4}
&&u_4=\frac{1}{4}\frac{-x+4a\sqrt t\tanh\left(y-\frac{x}{\sqrt t}\right)}{at},\\
\nonumber&&v_4=f_5(y).
\end{eqnarray}


\subsection{Reduction using $<a_4X_4+X_6(F_1)>$}
The basis associated with this subalgebra is $-2a_4t\frac{\partial}{\partial t}-a_4x\frac{\partial}{\partial x}-F_1(y)\frac{\partial}{\partial y}+a_4u\frac{\partial}{\partial u}+vF_1'(y)\frac{\partial}{\partial v}$ and the corresponding similarity variables are $m=\frac{t}{x^2}$ and $n=-\ln(x)+a_4\int{\frac{1}{F_1(y)}}dy.$ With the help of these similarity variables
we have the ansatz for $u$ and $v$ as
$$u(t,x,y)=\frac{U(m,n)}{x},v(t,x,y)=\frac{V(m,n)}{F_1(y)}.$$
Using this ansatz for $u,v$ and substituting in the given system \eqref{gov} we have the following reduced PDE system
\begin{eqnarray}\label{red3}
\nonumber&&a_4U_{mn}+4aa_4mU_mU_n+2aa_4U_n^2+4aa_4mUU_{mn}+2aa_4UU_{nn}+4aa_4UU_n+4aa_4m^2U_{mmn}+4aa_4mU_{mnn}\\
\nonumber&&+10aa_4mU_{mn}+aa_4U_{nnn}+3aa_4U_{nn}+2aa_4U_n+8bm^3V_{mmm}+12bm^2V_{mmn}+36m^2V_{mm}+6bmV_{mnn}\\
&&+24bmV_{mn}+24bmV_m+bV_{nnn}+3bV_{nn}+2bV_n=0,\\
\nonumber&&V_m-4cm^2V_{mm}-4cmV_{mn}-6cmV_m-cV_{nn}-cV_n+2gmUV_m+gUV_n=0.
\end{eqnarray}
In order to find the solution of the  reduced system \eqref{red3}, we again apply the Lie symmetry analysis and derive the ansatz as $U(m,n)=f_6(m),V(m,n)=\frac{n}{C_9}+f_7(m)$ where $C_9$ is an arbitrary constant and $f_6(m),f_7(m)$ are unknown functions those have to be determined. By exploiting $U$ and $V$ into the reduced system \eqref{red3} we have the following ODE system
\begin{eqnarray*}
&&4C_9m^3f_7'''+18C_9m^2f_7''+12C_9mf_7+1=0,\\
&&-C_9f_7'+4cC_9m^2f_7''+6cC_9mf_7'-2gC_9mf_6f_7'-gf_6+c=0.
\end{eqnarray*}

\begin{figure}[h!]
  \centering
  \subfloat[]{\includegraphics[scale=0.28]{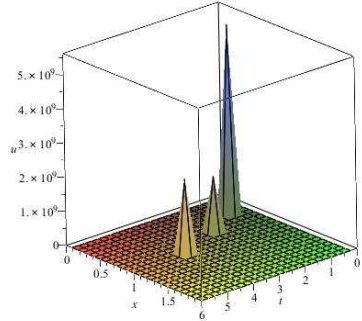}\label{sol5u2}}

  \subfloat[]{\includegraphics[scale=0.28]{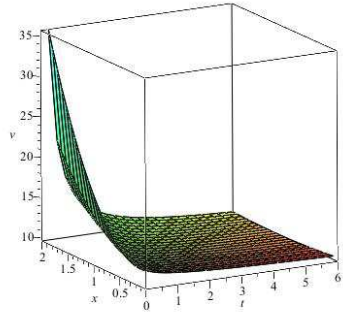}\label{sol5v2}}
  \subfloat[]{\includegraphics[scale=0.22]{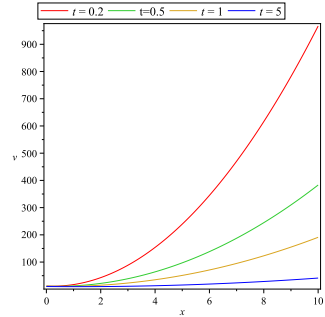}\label{sol5v22}}
  \caption{Solution profile of \eqref{gov} for a solution \eqref{sol5}: (a) 3d profile of $u_5$ (b) 3d profile of $v_5$ when $y=10$ (c) 2d profile of $v_5$ when $y=10$ for various values of $t$. }
  \label{fsol5}
  \end{figure}

By solving the above system of ODEs yields
$$f_6(m)=-\frac{1}{4}\frac{-2C_9C_{11}\sqrt m+4C_9C_{10}+(1+8cC_9C_{10})m}{C_9g\sqrt m(2C_{10}\sqrt m-C_{11}m)},f_7(m)=-\frac{2C_{11}}{\sqrt m}+\frac{2C_{10}}{m}-\frac{\ln(m)}{2C_9}+C_{12},$$
where $C_{10},C_{11}$ and $C_{12}$ are arbitrary constants.
Thus we have the exact solution of the given system \eqref{gov} as
\begin{eqnarray}\label{sol5}
&&u_5=-\frac{1}{4}\frac{x(-2C_9C_{11}x\sqrt t+4C_9C_{10}x^2+t+8cC_9C_{10}t)}{C_9gt(-C_{11}x\sqrt t+2C_{10}x^2)},\\
\nonumber&&v_5=-\frac{1}{2}\frac{2t\ln(x)-2ta_4\int{\left(\frac{1}{F_1(y)}\right)}dy+4C_9C_{11}x\sqrt t-4C_{10}x^2+t\ln\left(\frac{t}{x^2}\right)-2C_9C_{12}t}{C_9tF_1(y)}.
\end{eqnarray}
We now discuss the physical significance of the solution profile \eqref{sol5} (see, Figure \ref{fsol5}) which represents $u_5$ and $v_5$ by considering $c=1,g=1,a_4=1,C_9=1,C_{10}=1,C_{11}=1,C_{12}=1$ and by choosing $F_1(y)=1$. The solution profile of $u_5$ is illustrated in the Figure \ref{sol5u2} which represents a 3-lump type soliton. The 3-dimensional profile of $v_5$ is depicted in the Figure \ref{sol5v2} at fixed $y=10$ and the corresponding 2-dimensional profile is shown in the Figure \ref{sol5v22} with respect to $x$ for different values of $t$. Here Figure \ref{sol5v2} shows that $v_5$ increases rapidly if $x$ increases at initial time period and then it decreases gradually as time evolves. From the corresponding 2-dimensional profile \ref{sol5v22} we observed that $v_5$ increases with respect to $x$ for any $t$ but as $t$ increases, the rate of increasing  of $v_5$ gradually decreases and after some time (say, t=5) the curve becomes almost horizontal line.


Another solution of the reduced PDE system \eqref{red3} is $U(m,n)=f_8(m),V(m,n)=C_{13}$ where $C_{13}$ is an arbitrary constant and $f_8(m)$ is an arbitrary function. Thus we have another solution of \eqref{gov} as follows:
\begin{eqnarray}\label{sol6}
&&u_6=f_8\left(\frac{t}{x^2}\right),v_6=\frac{C_{13}}{F_1(y)}.
\end{eqnarray}
We obtain another two solutions of the reduced system \eqref{red3} as
\begin{eqnarray*}
&&U(m,n)=-\frac{1}{2}\frac{2C_{15}+C_{16}\sqrt m+4cC_{15}m}{g(2C_{15}\sqrt m+C_{16}m)},\\
&&V(m,n)=C_{14}+\frac{C_{15}}{m}+\frac{C_{16}}{\sqrt m}
\end{eqnarray*}
and
\begin{eqnarray*}
&&U(m,n)=\frac{c(-2+n)}{gn},\\
&&V(m,n)=\frac{1}{n},~~~\text{if}~~a_4=\frac{1}{2}\frac{bg^2}{ac(2c+g)}
\end{eqnarray*}
where $C_{14},C_{15}$ and $C_{16}$ are arbitrary constants.

\begin{figure}[h!]
  \centering
  \subfloat[]{\includegraphics[scale=0.28]{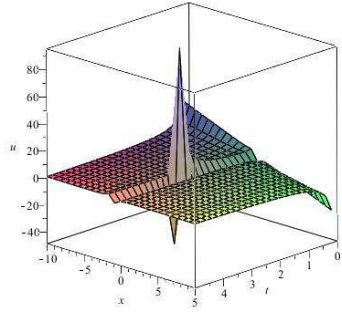}\label{sol7u2}}

  \subfloat[]{\includegraphics[scale=0.28]{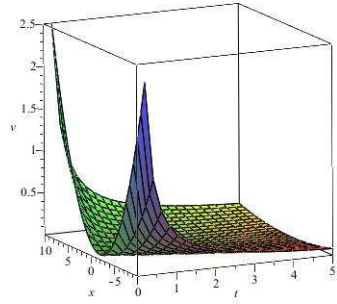}\label{sol7v2}}
  \subfloat[]{\includegraphics[scale=0.22]{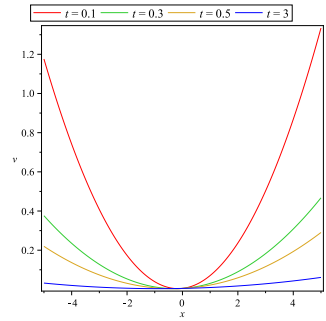}\label{sol7v22}}
  \caption{Solution profile of \eqref{gov} for a solution \eqref{sol7}: (a) 3d profile of $u_7$ (b) 3d profile of $v_7$ when $y=10$ (c) 2d profile of $v_7$ when $y=10$ for various values of $t$. }
  \label{fsol7}
  \end{figure}

Using the above two solutions of the reduced PDE system \eqref{red3} for $U(m,n)$ and $V(m,n)$ we obtain two more exact solutions for the given system \eqref{gov} as follows:

\begin{eqnarray}\label{sol7}
&&u_7=-\frac{x}{2t}\frac{2C_{15}x^2+C_{16}x\sqrt t+4cC_{15}t}{g(2C_{15}x^2+C_{16}x\sqrt t)},\\
\nonumber&&v_7=\frac{C_{14}t+C_{15}x^2+C_{16}x\sqrt t}{tF_1(y)}
\end{eqnarray}
\begin{figure}[h!]
  \centering
  \subfloat[]{\includegraphics[scale=0.28]{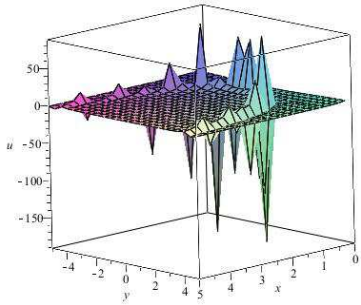}\label{sol8u}}
  \hfill
  \subfloat[]{\includegraphics[scale=0.28]{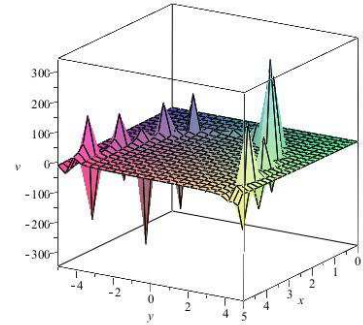}\label{sol8v}}
  \caption{Solution profile of \eqref{gov} for a solution \eqref{sol8}:  (a) 3d profile of $u_8$  (b) 3d profile of $v_8$}
  \label{fsol8}
  \end{figure}


and
\begin{eqnarray}\label{sol8}
&&u_8=\frac{c\left(2+\ln(x)-a_4\left(\int{\frac{1}{F_1(y)}}dy\right)\right)}{gx\left(\ln(x)-a_4\left(\int{\frac{1}{F_1(y)}}dy\right)\right)},\\
\nonumber&&v_8=-\frac{1}{\left(\ln(x)-a_4\left(\int{\frac{1}{F_1(y)}}dy\right)\right)F_1(y)},~~\text{if}~~a_4=\frac{1}{2}\frac{bg^2}{ac(2c+g)}.
\end{eqnarray}

Now we discuss the physical behavior of the solution profile  \eqref{sol7} in the Figure \ref{fsol7} for $u_7$ and $v_7$ by choosing the parameters $c=1,g=1,C_{14}=1,C_{15}=1,C_{16}=1$ and considering $F_1(y)=2y^2$. The Figure \ref{sol7u2} indicates a multiple breather soliton type solution for $u_7$. The 3-dimensional profile of $v_7$ is depicted in the Figure \ref{sol7v2} by choosing $y = 10$ and Figure \ref{sol7v22} illustrates the 2-dimensional profile of $v_7$ at fixed $y = 10$ with respect to $x$ for varying $t$ which indicates an upward parabola.


The physical significance of the solution profile \eqref{sol8} is depicted in the Figure \ref{fsol8} by considering the parameters $a=1,b=1,c=1,g=1,a_4=\frac{1}{2}\frac{bg^2}{ac(2c+g)}=\frac{1}{6}$ and by setting $F_1(y)=\frac{1}{y}$. The surface profiles of $u_8$ and $v_8$ are illustrated in the Figure \ref{sol8u} and Figure \ref{sol8v} respectively those indicate a multiple breather soliton type solutions.

\subsection{Reduction using $<X_6(F_1)>$}
The representative for this class of subalgebra is given by $X_6(F_1)=-F_1(y)\frac{\partial}{\partial y}+vF_1'(y)\frac{\partial}{\partial v}$. The governing system \eqref{gov} with the similarity transformations $u(t,x,y)=U(t,x),v(t,x,y)=\frac{V(t,x)}{F_1(y)}$
%

\begin{figure}[h!]
  \centering
  \subfloat[]{\includegraphics[scale=0.30]{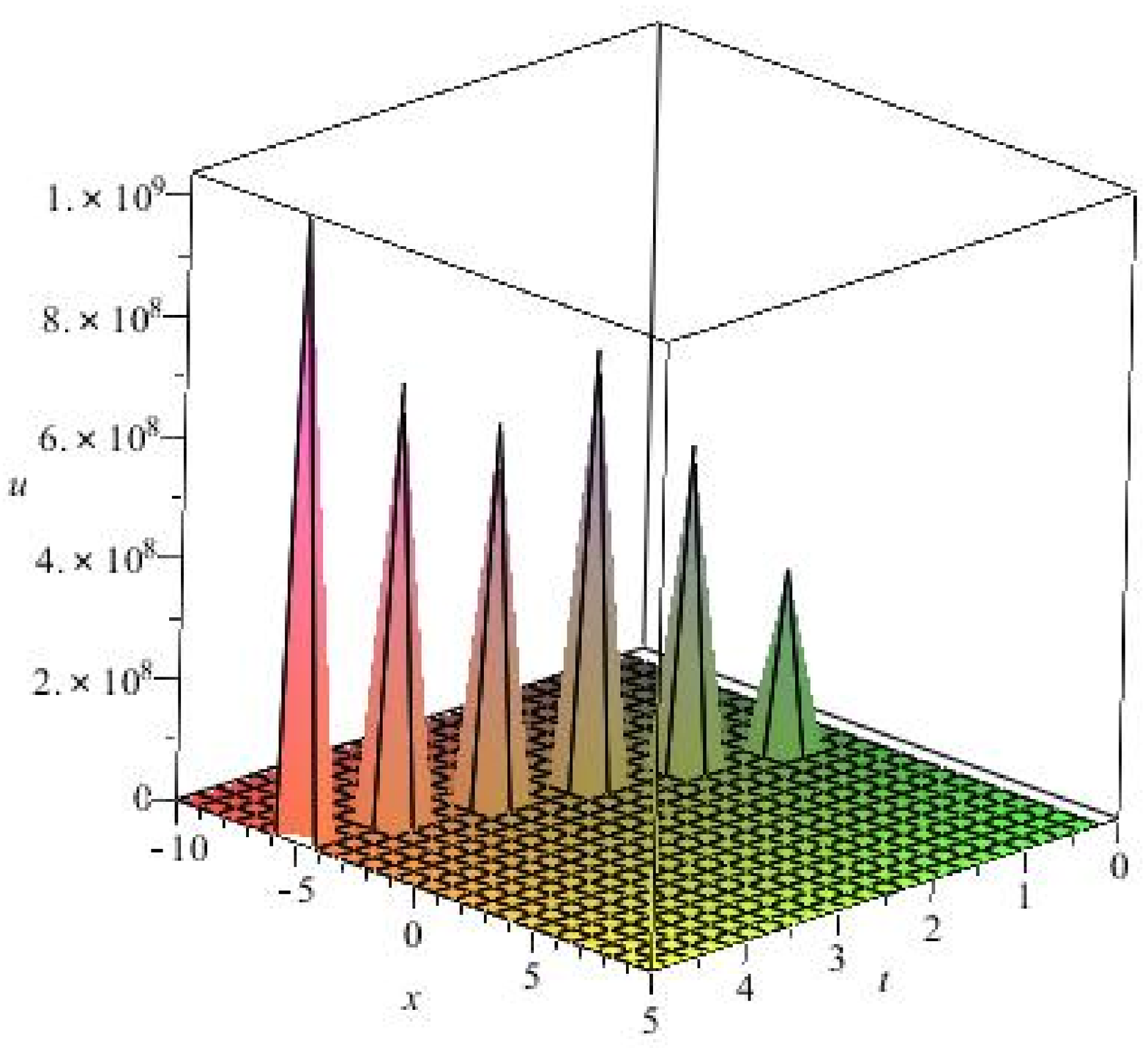}\label{sol9u2}}

  \subfloat[]{\includegraphics[scale=0.28]{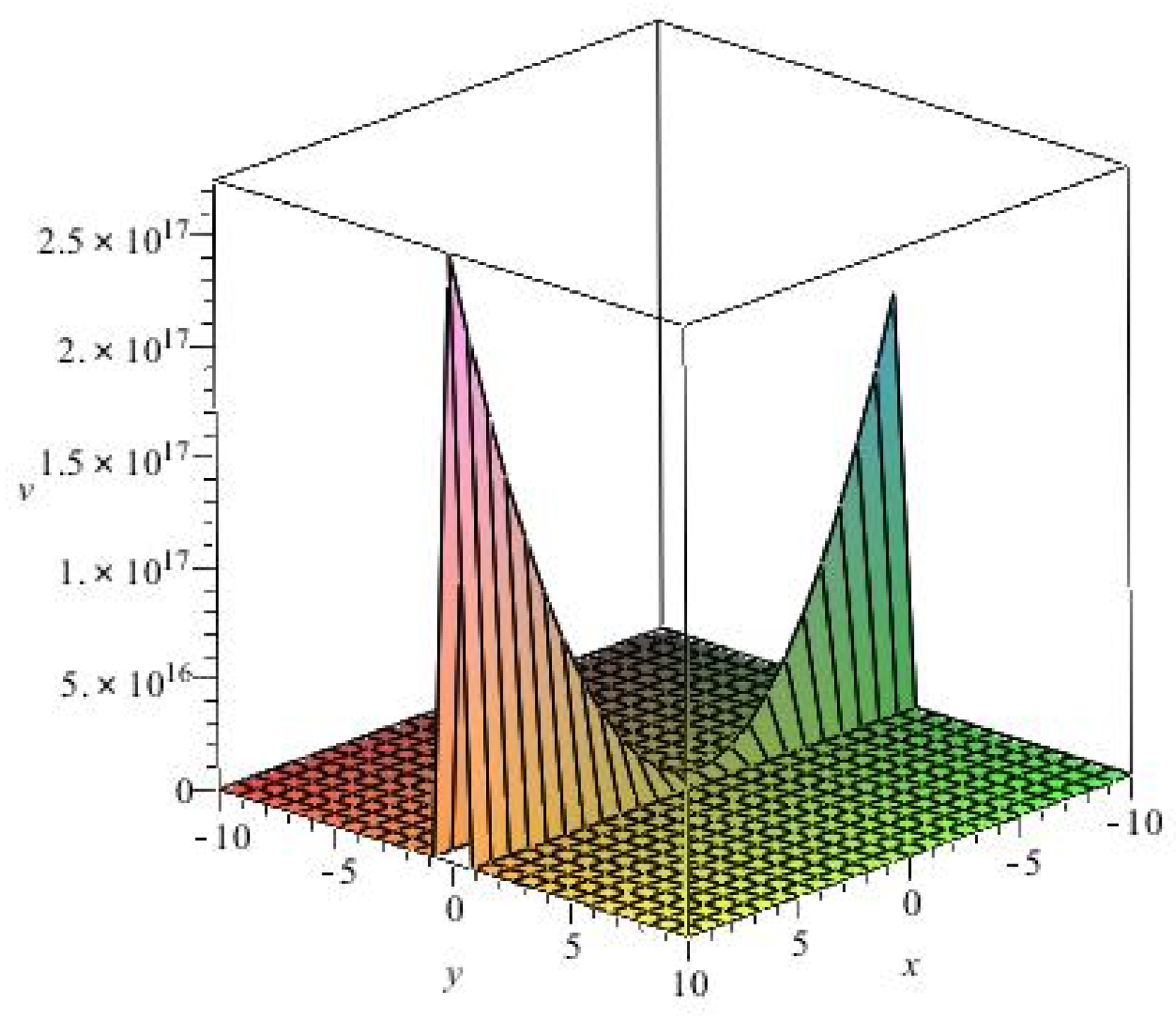}\label{sol9v1}}
  \subfloat[]{\includegraphics[scale=0.28]{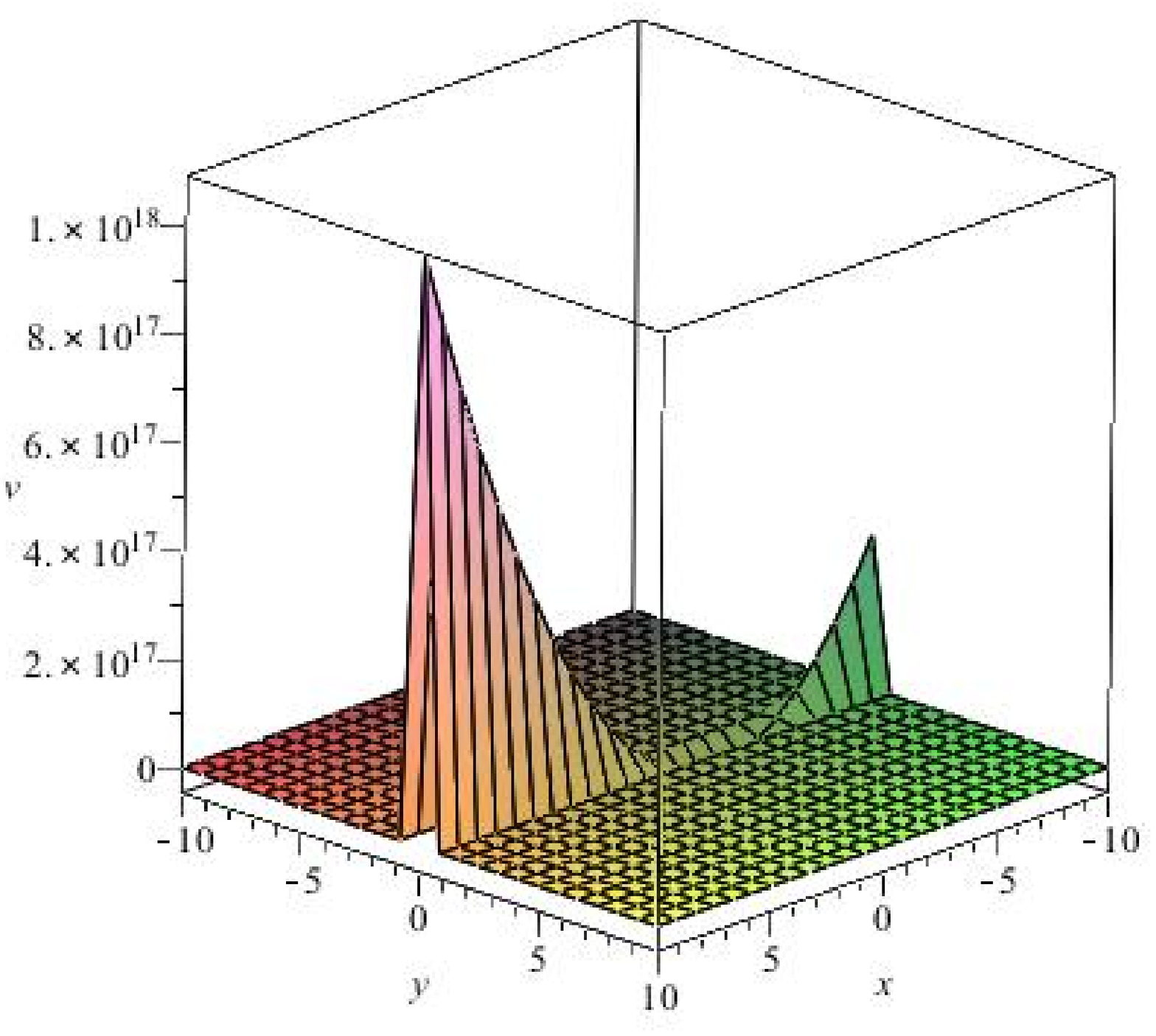}\label{sol9v2}}
  \subfloat[]{\includegraphics[scale=0.28]{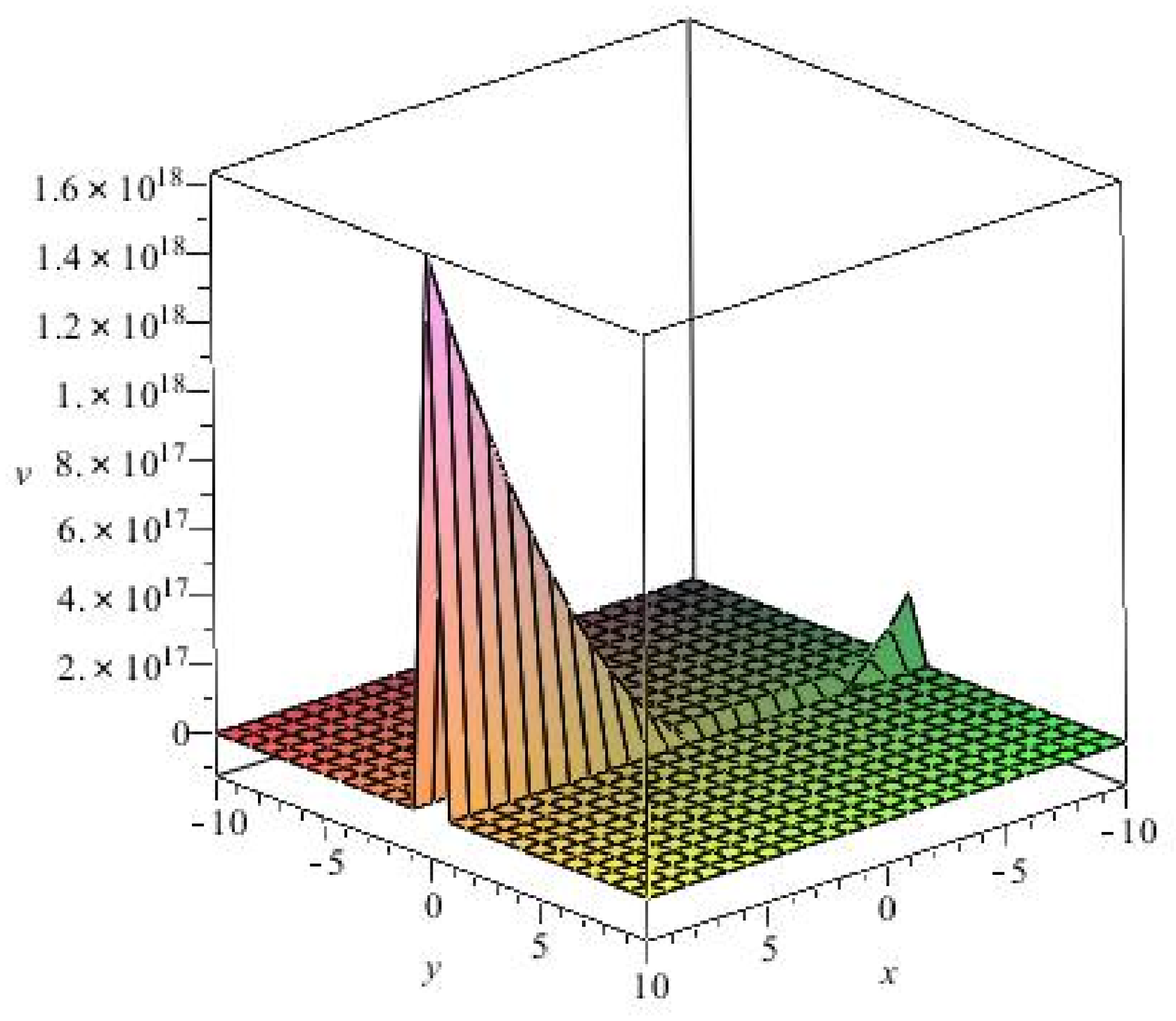}\label{sol9v3}}

  \subfloat[]{\includegraphics[scale=0.28]{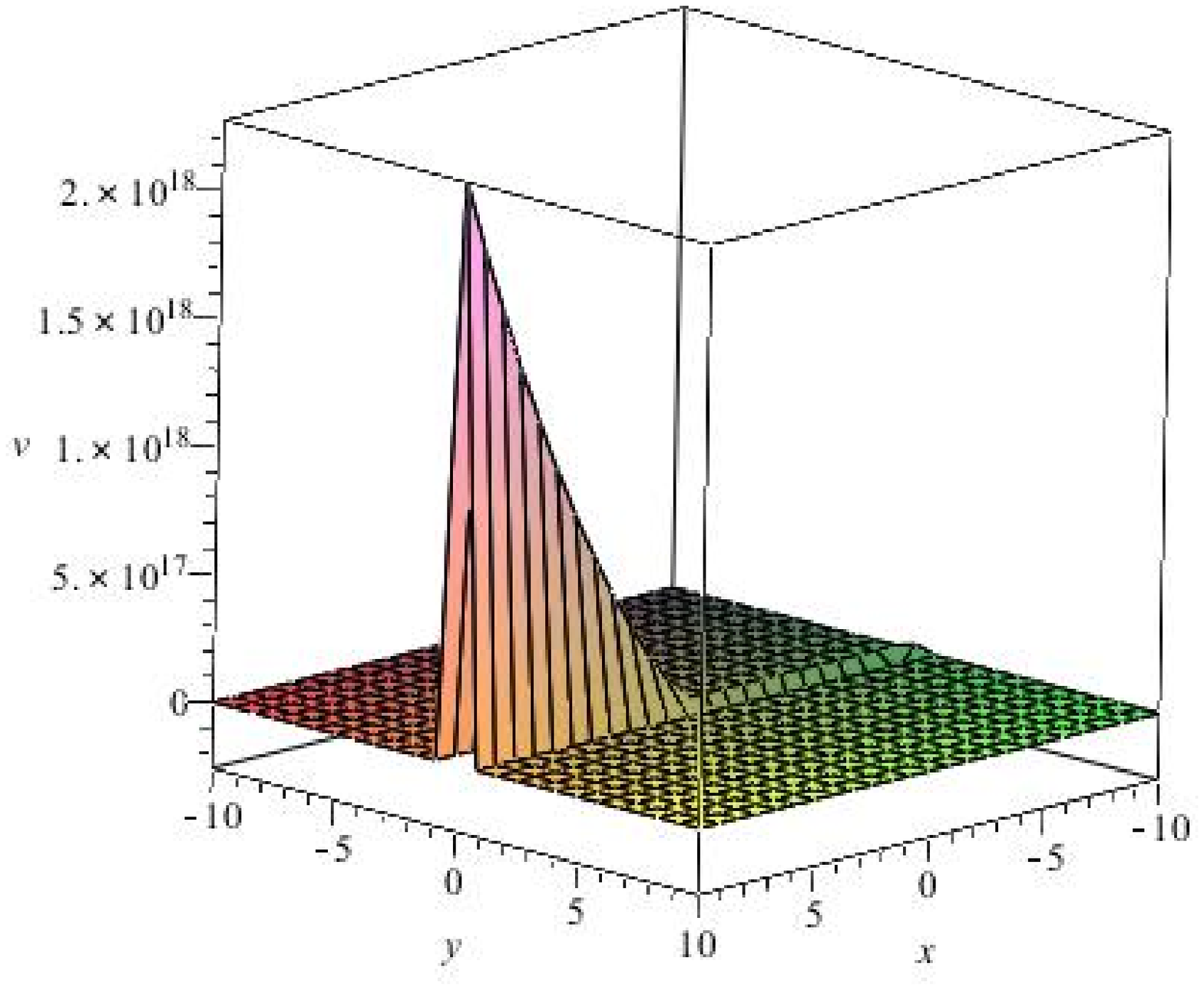}\label{sol9v4}}
  \subfloat[]{\includegraphics[scale=0.28]{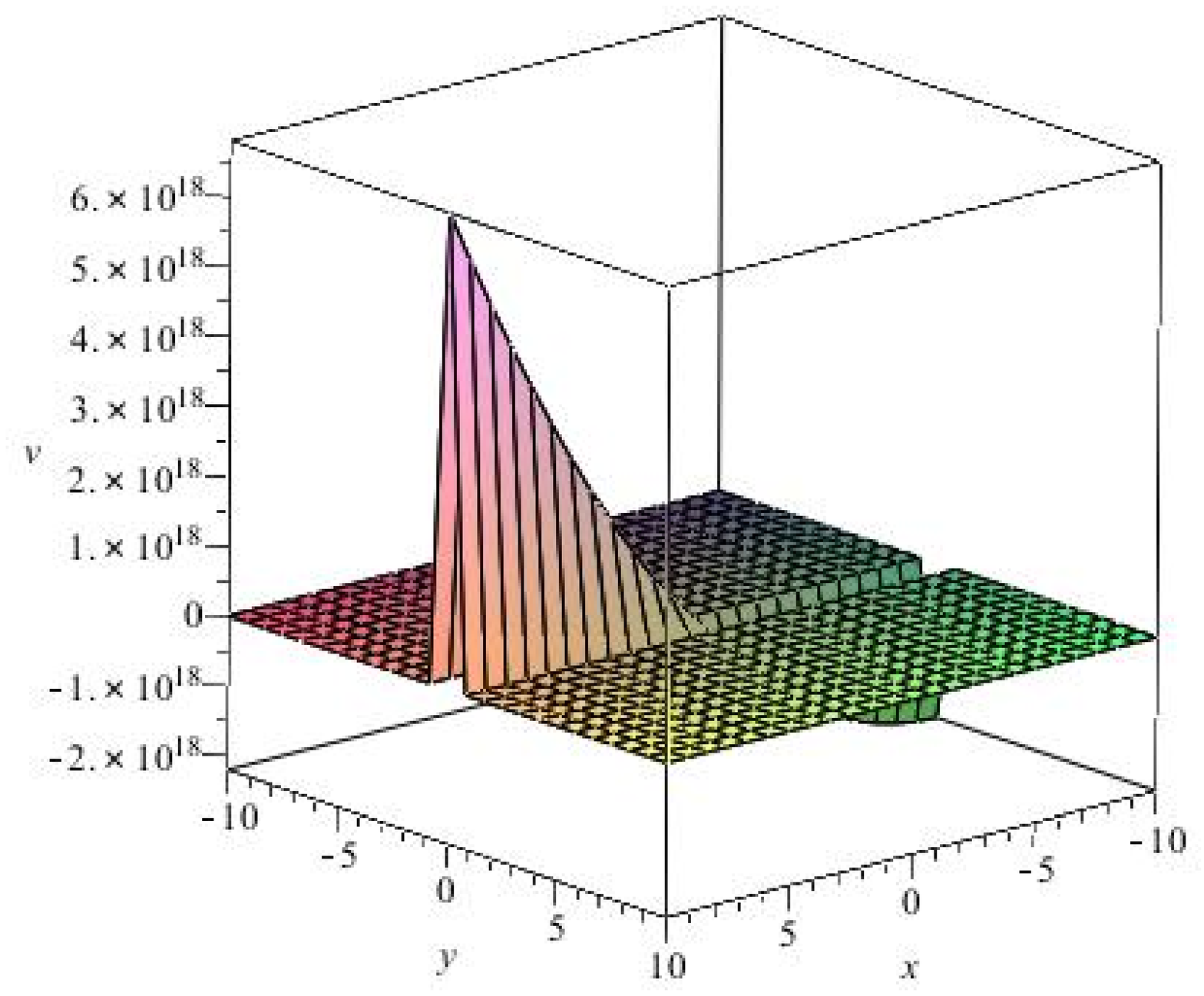}\label{sol9v5}}
  \subfloat[]{\includegraphics[scale=0.28]{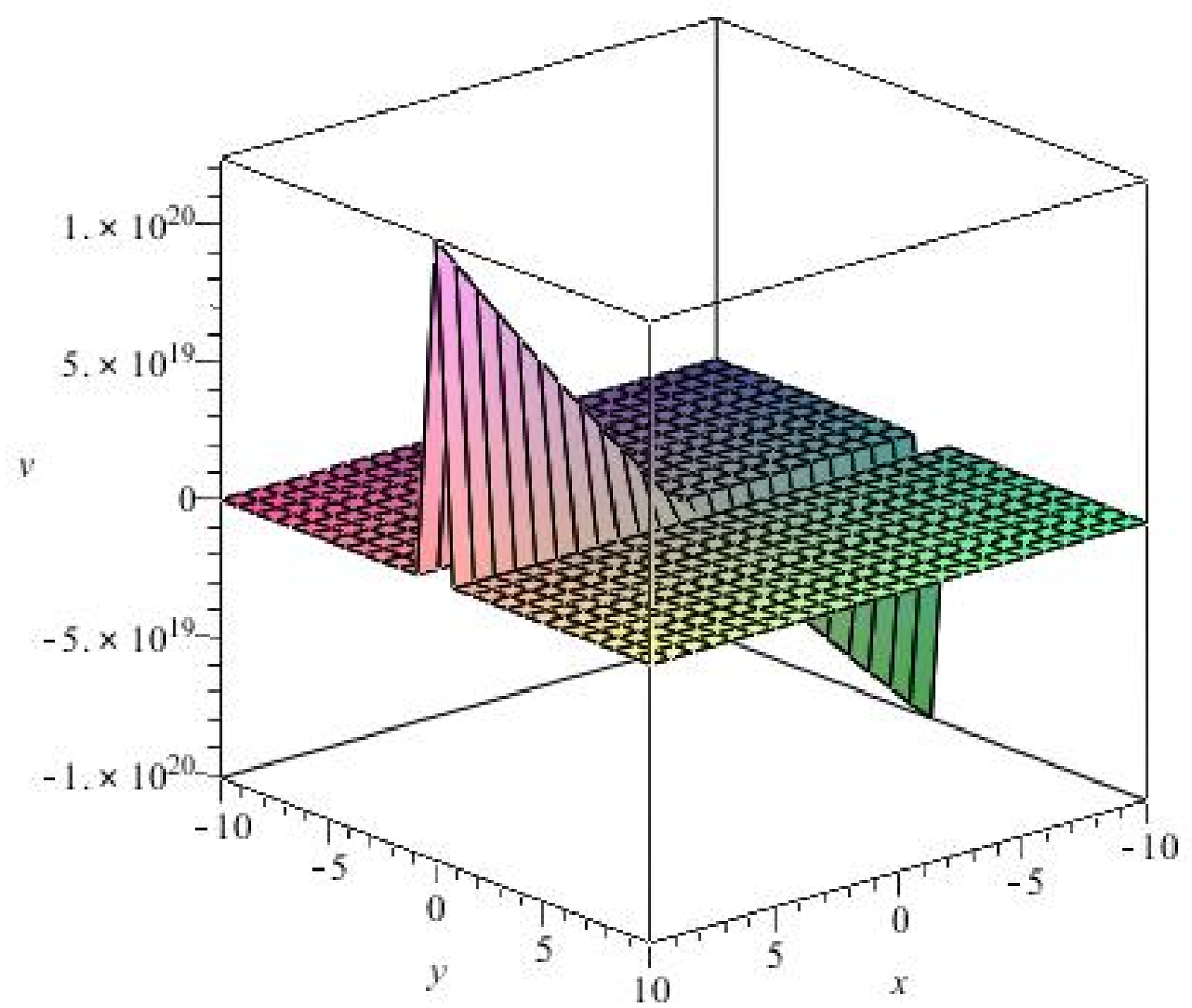}\label{sol9v6}}
  \caption{Solution profile of \eqref{gov} for a solution \eqref{sol9}: (a) $u_{9}$ , (b) $v_{9}$ at $t=1$, (c) $v_{9}$ at $t=3$, (d) $v_{9}$ at $t=4$, (e) $v_{9}$ at $t=5$, (f) $v_{9}$ at $t=10$, (g) $v_{9}$ at $t=50$}
  \label{fsol9}
  \end{figure}

reduces to the following  PDE system
\begin{eqnarray}\label{red4}
&&V_{xxx}=0,\\
\nonumber&&V_t-cV_{xx}-gUV_x=0.
\end{eqnarray}
We easily solve this system \eqref{red4} and obtain
 $$U(t,x)=\frac{1}{2}\frac{f_9'(t)x^2+2f_{10}'(t)x+2f_{11}'(t)-2cf_9(t)}{g(f_9(t)x+f_{10}(t))},V(t,x)=\frac{1}{2}f_9'(t)x^2+f_{10}(t)x+f_{11}(t)$$
where $f_9,f_{10}$ and $f_{11}$ are arbitrary functions of $t$. Thus finally we have the solution for the given system of PDEs of the form
\begin{eqnarray}\label{sol9}
&&u_9=\frac{1}{2}\frac{f_9'(t)x^2+2f_{10}'(t)x+2f_{11}'(t)-2cf_9(t)}{g(f_9(t)x+f_{10}(t))},\\
\nonumber&&v_9=\frac{\frac{1}{2}f_9'(t)x^2+f_{10}(t)x+f_{11}(t)}{F_1(y)}.
\end{eqnarray}

Now by choosing the parameters $c=1,g=1$ and  considering $F_1(y)=y,f_9(t)=t,f_{10}(t)=t^2$ and $f_{11}(t)=t$ we study the physical behavior of $u_9$ and $v_9$ given in \eqref{sol9} (see, Figure \ref{fsol9}). The Figure  \ref{sol9u2} demonstrates a 5-lump type soliton or multiple lump type soliton profile for $u_9$. On the other hand we draw the 3-dimensional  surface for $v_9$ with respect to the spatial variables $x$ and $y$ for different values of $t$ in the Figure \ref{sol9v1}-\ref{sol9v6}. It is very interesting to observe that initially at $t=1$, the corresponding Figure \ref{sol9v1} indicates 2-soliton profile where both the peaks are in the same (positive) direction. After some time, say $t=5$, it annihilates into a 1-soliton or single soliton (see, Figure \ref{sol9v4}) and further later, say at $t=50$, the profile $v_9$ again behaves like a 2-soliton solution (see, Figure \ref{sol9v6}) but in the opposite direction.

\section{Traveling wave solutions}\label{s5}
The exploration of the traveling wave solutions, in particular soliton solutions of nonlinear system of PDEs play a vital role in describing the characters of nonlinear problems in the area of engineering, applied science and mathematical physics. It also describes many interesting physical phenomena in the study of dynamical systems. Here we consider the traveling wave solution of the form $$u(t,x,y)=U\left(\frac{l_1x-l_2t}{l_1},\frac{l_1y-l_3t}{l_1}\right),v(t,x,y)=V\left(\frac{l_1x-l_2t}{l_1},\frac{l_1y-l_3t}{l_1}\right)$$
which is invariant under the symmetry $l_2X_1+l_3X_2+l_1X_3.$ Using this form of $u,v$ and exploiting them into the given system \eqref{gov} which yields the reduced system of PDEs
\begin{eqnarray}\label{red5}
&&l_2U_{mn}+l_3U_{nn}+2al_1U_mU_n+2al_1UU_{mn}-al_1U_{mmn}+bl_1V_{mmm}=0,\\
\nonumber&&l_2V_m+l_3V_n+cl_1V_{mm}+gl_1UV_m=0
\end{eqnarray}
where $m=\frac{l_1x-l_2t}{l_1}$ and $n=\frac{l_1y-l_3t}{l_1}.$
In order to solve the system \eqref{red5}, we use the Lie symmetry approach and obtain the solutions of \eqref{red5} in the form of $U(m,n)=f_{12}(m),V(m,n)=\frac{C_{19}n}{C_{18}}+f_{13}(m)$ and $U(m,n)=f_{14}(n),V(m,n)=\frac{C_{19}m}{C_{17}}+f_{15}(n)$ where $C_{17},C_{18}$ and $C_{19}$ are arbitrary constants and $f_{12}(m),f_{13}(m),f_{14}(n)$ and $f_{15}(n)$ are unknown functions.

\begin{figure}[h!]
  \centering
  \subfloat[]{\includegraphics[scale=0.28]{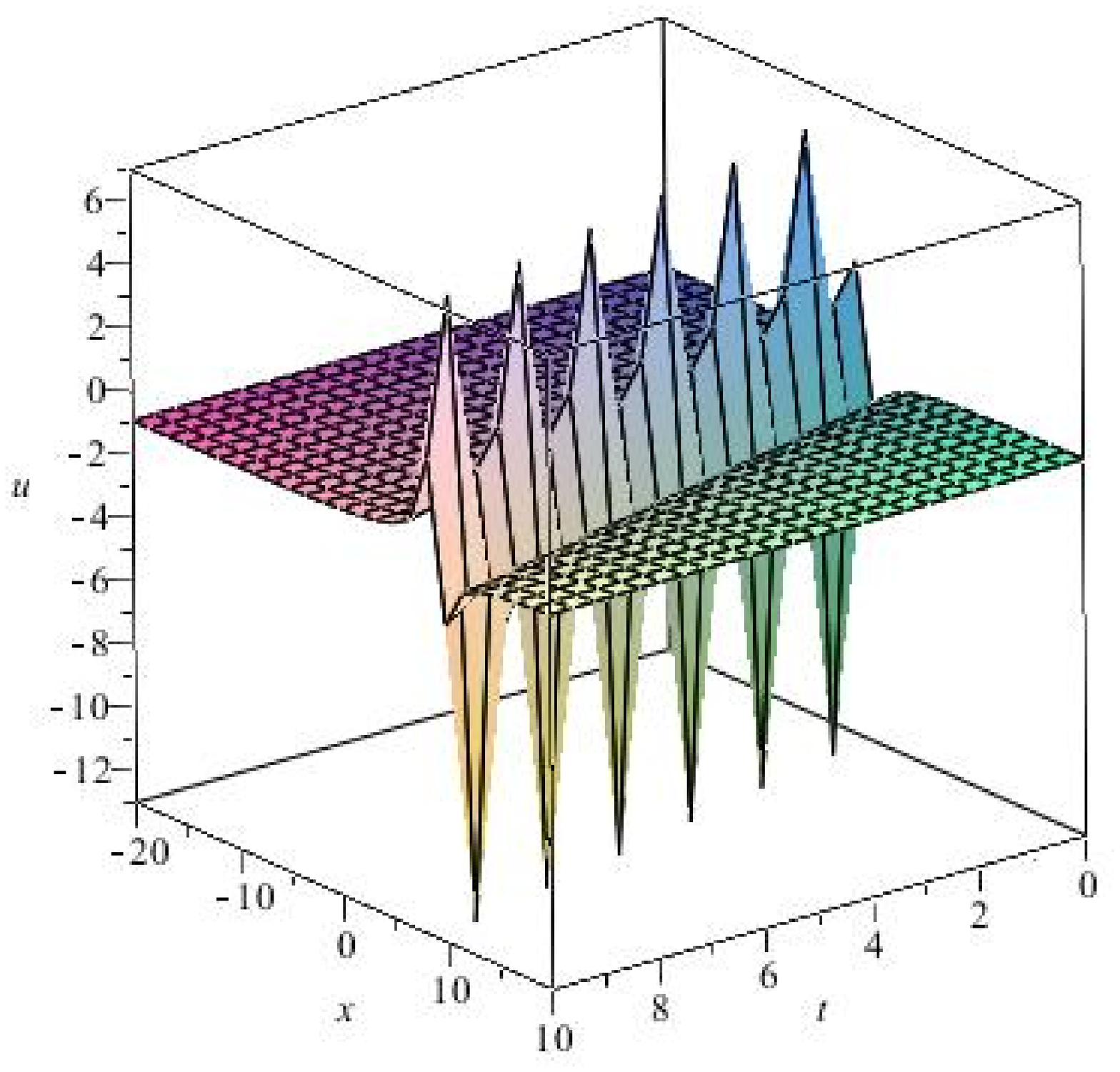}\label{sol10u2}}
   \subfloat[]{\includegraphics[scale=0.28]{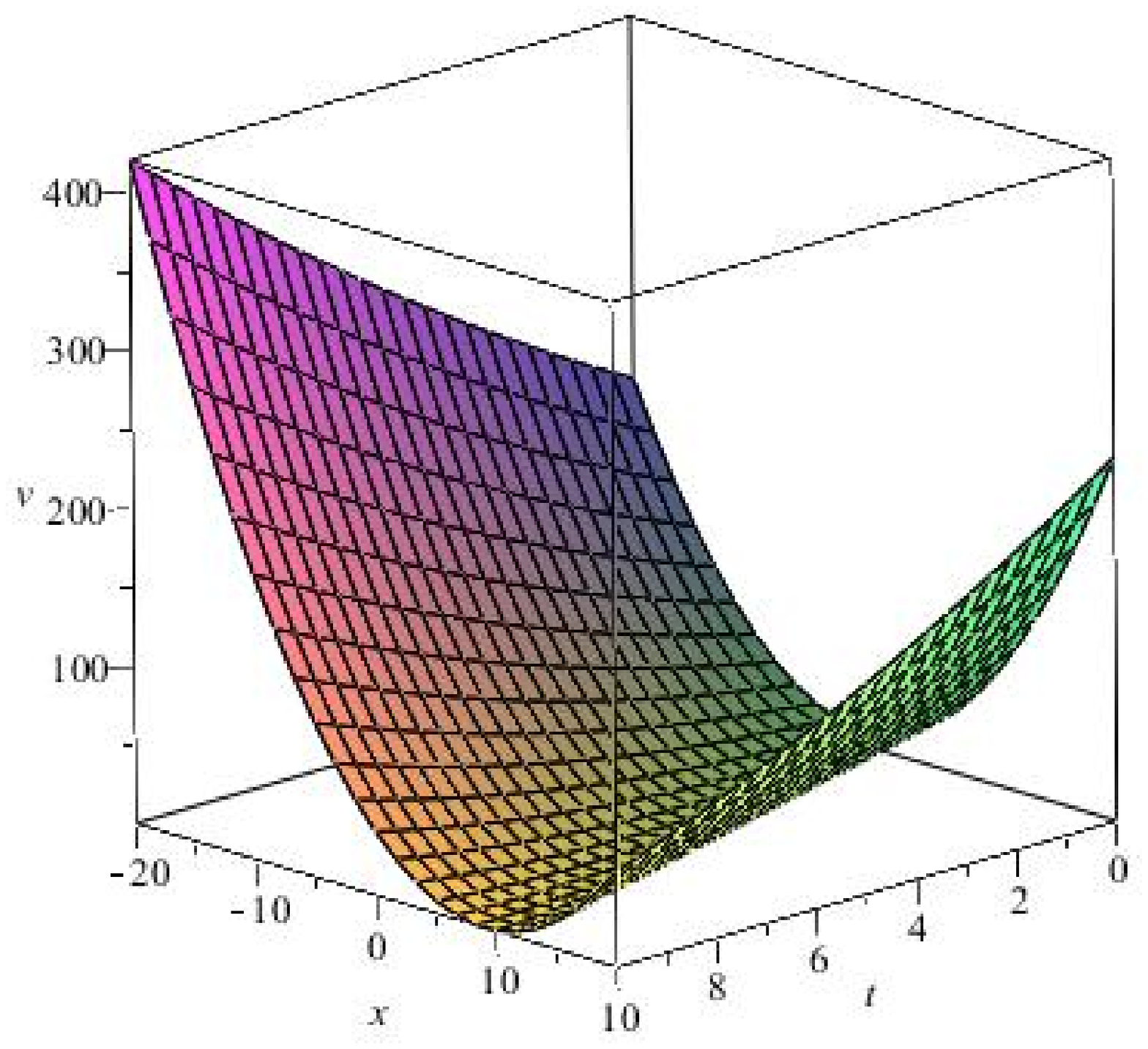}\label{sol10v2}}
  \caption{Solution profile of \eqref{gov} for a solution \eqref{sol10}: (a) 3d profile of $u_{10}$  (b) 3d profile of $v_{10}$ when $y=10$}
  \label{fsol10}
  \end{figure}

 Then substituting these forms of $U$ and $V$ into the reduced system \eqref{red5} and solving the corresponding ODE systems we obtain
$$f_{12}(m)=-\frac{l_2C_{18}C_{20}m+l_2C_{18}C_{21}+l_3C_{19}+cl_1C_{18}C_{20}}{gl_1C_{18}(C_{20}m+C_{21})},~~f_{13}(m)=\frac{1}{2}C_{20}m^2+C_{21}m+C_{22}$$
and
$$f_{14}(n)=C_{24}n+C_{25},~~f_{15}(n)=\frac{1}{2}\frac{-gl_1C_{19}C_{24}n^2-2C_{19}n(l_2+gl_1C_{25})+2l_3C_{17}C_{23}}{l_3C_{17}}$$
where $C_{20},...,C_{25}$ are integration constants which in turn yields the following solutions of the given system
\begin{eqnarray}\label{sol10}
&&u_{10}=\frac{l_2C_{18}C_{20}(l_1x-l_2t)+l_1l_2C_{18}C_{21}+l_1l_3C_{19}+cl_1^2C_{18}C_{20}}{l_1gC_{18}\{C_{20}(l_2t-l_1x)-l_1C_{21}\}},\\
\nonumber&&v_{10}=\frac{1}{2}\frac{2l_1C_{19}(l_1y-l_3t)+C_{18}C_{20}(l_2^2t^2+l_1^2x^2)-2l_1l_2C_{18}C_{20}tx+2C_{18}C_{21}l_1(l_1x-l_2t)
+2l_1^2C_{18}C_{22}}{l_1^2C_{18}}
\end{eqnarray}
and
\begin{eqnarray}\label{sol11}
&&u_{11}=\frac{C_{24}(l_1y-l_3t)+l_1C_{25}}{l_1},\\
\nonumber&&v_{11}=\frac{1}{2}\frac{2l_1C_{19}(l_3x-l_2y)-gC_{19}C_{24}(l_3^2t^2+l_1^2y^2)+2gl_1C_{19}C_{25}(l_3t-l_1y)+2l_1l_3C_{17}C_{23}}{l_1l_3C_{17}}.
\end{eqnarray}
We now discuss the physical significance of the solution profile of \eqref{sol10} by choosing the parameters $c=1,g=1,l_1=1,l_2=1,l_3=1,C_{18}=1,C_{19}=1,C_{20}=1,C_{21}=1$ and $C_{22}=1$ in the Figure \ref{fsol10}. Figure \ref{sol10u2} demonstrates the singular kink traveling wave solution profile for $u_{10}$. Whilst we illustrate the surface profile of $v_{10}$ in the Figure \ref{sol10v2} by choosing $y=10$  which indicates an upward parabola.


We compute some other solutions of the reduced PDE system \eqref{red5} as given by
$$U(m,n)=\frac{2cC_{27}\tanh\left(C_{26}+C_{27}m-\frac{l_2C_{27}n}{l_3}\right)}{g},V(m,n)=C_{28}-\frac{2acl_2C_{27}(2c+g)\tanh\left(C_{26}+C_{27}m-\frac{l_2C_{27}n}{l_3}\right)}
{bg^2l_3}$$
and

\begin{figure}[h!]
  \centering
  \subfloat[]{\includegraphics[scale=0.30]{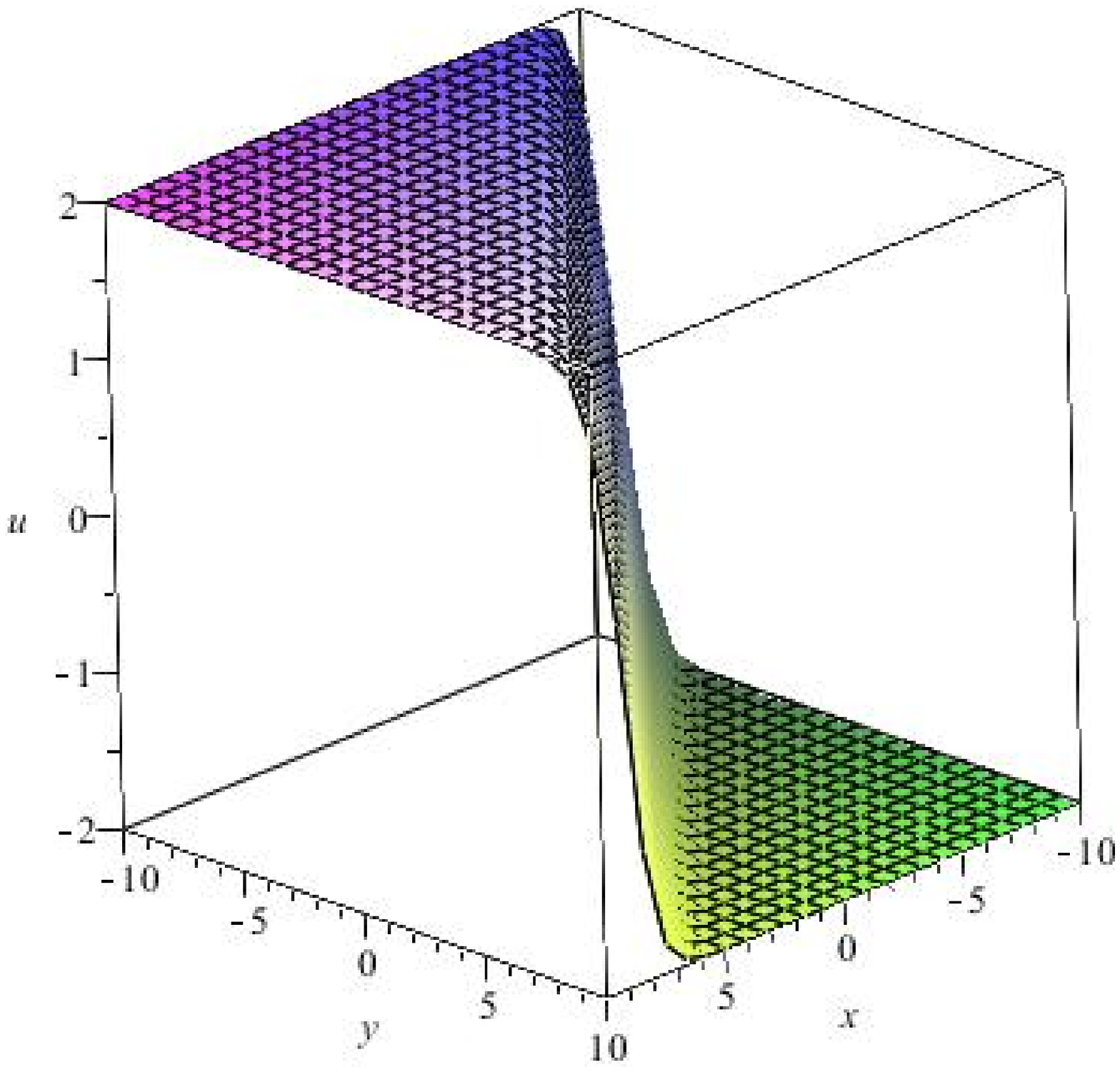}\label{sol12u}}
  \hfill
  \subfloat[]{\includegraphics[scale=0.30]{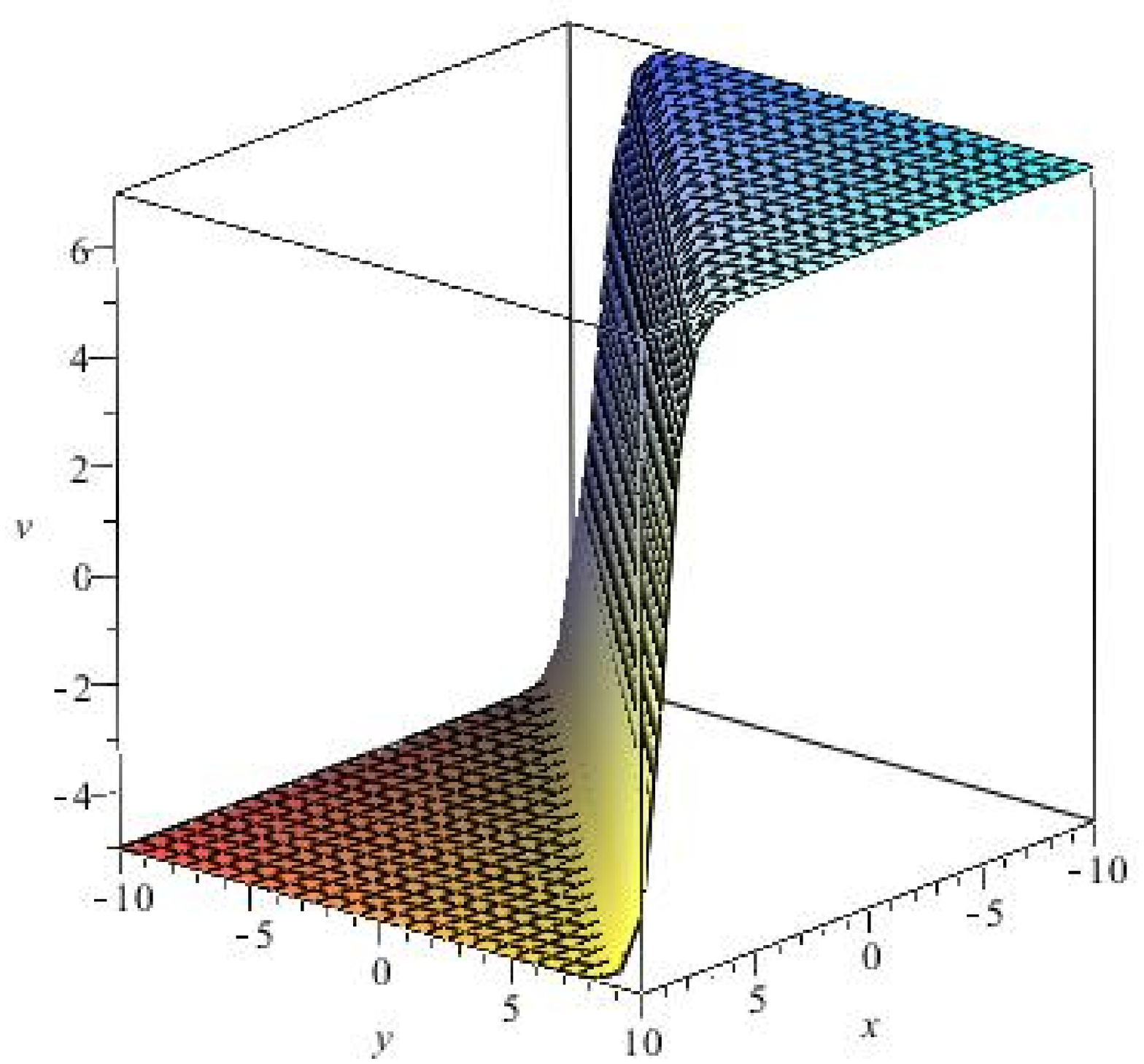}\label{sol12v}}
  \caption{Solution profile of \eqref{gov} for a solution \eqref{sol12}: (a) 3d profile of $u_{12}$  (b) 3d profile of $v_{12}$}
  \label{fsol12}
  \end{figure}

$$U(m,n)=-\frac{1}{2}\frac{l_2C_{31}+l_3C_{32}}{al_1C_{31}}-C_{31}\tanh(C_{31}m+C_{32}n+C_{30}),V(m,n)=C_{29}$$
where $C_{26},...,C_{32}$ are arbitrary constants. This again yields to new exact solutions of the given system \eqref{gov} as follows:
\begin{eqnarray}\label{sol12}
&&u_{12}=\frac{2cC_{27}\tanh\left(\frac{C_{27}(l_3x-l_2y)+l_3C_{26}}{l_3}\right)}{g},\\
\nonumber&&v_{12}=\frac{bg^2l_3C_{28}-2acC_{27}l_2(2c+g)\tanh\left(\frac{C_{27}(l_3x-l_2y)+l_3C_{26}}{l_3}\right)}{bg^2l_3}
\end{eqnarray}


and
\begin{eqnarray}\label{sol13}
&&u_{13}=-\frac{1}{2}\frac{l_2C_{31}+l_3C_{32}+2al_1C_{31}^2\tanh\left(\frac{l_1C_{30}-(l_2C_{31}+l_3C_{32})t+l_1C_{31}c+l_1C_{32}y}{l_1}\right)}
{al_1C_{31}},\\
\nonumber&&v_{13}=C_{29}.
\end{eqnarray}
We study the physical significance of the stationary solution profile given in \eqref{sol12} by considering the parameters $a=1,c=1,g=1,l_2=1,l_3=1,C_{26}=1,C_{27}=1,C_{28}=1$ and we demonstrate their 3-dimensional profiles in the Figure \ref{fsol12}. Here we observed that $u_{12}$ behaves like an anti-kink type soliton (see, Figure \ref{sol12u}) and $v_{12}$ satisfies the properties of kink type soliton profile (see, Figure \ref{sol12v}).

In \cite{wang20202+} authors obtained only one stationary domain walls solution of \eqref{gov} of the following form
\begin{eqnarray}\label{recover}
&&u(t,x,y)=A_1\tanh[B_1x+B_2y],\\
\nonumber&&v(t,x,y)=A_2\tanh[B_1x+B_2y]
\end{eqnarray}
where $bA_2B_1^2=aA_1(A_1+B_1)B_2$ and $gA_1=2cB_1$. One can observe that the solution \eqref{recover} of \eqref{gov} is a particular case of our solution \eqref{sol12} by considering $C_{26}=0, C_{28}=0$ in \eqref{sol12} and letting $A_1=\frac{2cC_{27}}{g},A_2=-\frac{2acC_{27}l_2(2c+g)}{bg^2l_3},B_1=C_{27},B_2=-\frac{C_{27}l_2}{l_3}$ in \eqref{recover}.

 \par

Now we consider the solution of \eqref{gov} of the form $u(t,x,y)=U(k_1x+k_2y+k_3t+k_4),v(t,x,y)=V(k_1x+k_2y+k_3t+k_4)$ which yields a reduced system of ODEs of the form
\begin{eqnarray}\label{red6}
&&k_2k_3U''(h)-2al_1l_2U'(h)^2-2al_1l_2UU''(h)+al_1^2l_2U'''(h)-bl_1^3V'''(h)=0,\\
\nonumber&&l_3V'(h)-cl_1^2V''(h)-l_1gUV'=0
\end{eqnarray}
where $h=k_1x+k_2y+k_3t+k_4$ and $k_1,...,k_4$ are arbitrary constants. In general it is difficult to solve the above ODE system \eqref{red6} but we can obtain some particular class of exact solutions. One of them is  as follows:
\begin{eqnarray*}
&&U(h)=\frac{k_3C_{33}h+k_3C_{34}gk_1^2C_{33}}{k_1g(C_{33}h+C_{34})},\\
&&V(h)=\ln(C_{33}h+C_{34})
\end{eqnarray*}
provided some restrictions on the parameters involved in the given system \eqref{gov} such as $c=-g$ and $a=\frac{1}{2}\frac{g(bk_1^2+k_2k_3)}{k_2k_3}$. In this context, we obtain the following exact solution of
\begin{eqnarray}\label{sol14}
&&u_{14}=\frac{C_{33}k_3(k_1x+k_2y+k_3t+k_4)+k_3C_{34}-gk_1^2C_{33}}{k_1g(C_{33}k_1x+C_{33}k_2y+C_{33}k_3t+k_4+C_{34})},\\
\nonumber&&v_{14}=\ln(C_{33}(k_1x+k_2y+k_3t+k_4)+C_{34})
\end{eqnarray}


where $C_{33}$ and $C_{34}$ are arbitrary constants.

On the other hand we have another solution of the reduced system of ODEs \eqref{red6} as
\begin{eqnarray*}
&&U(h)=\frac{k_3\cosh(C_{35}h+C_{36})+2ck_1^2C_{35}\sinh(C_{35}h+C_{36})}{gk_1\cosh(C_{35}h+C_{36})},\\
&&V(h)=\tanh(C_{35}h+C_{36})
\end{eqnarray*}

\begin{figure}[h!]
  \centering
  \subfloat[]{\includegraphics[scale=0.30]{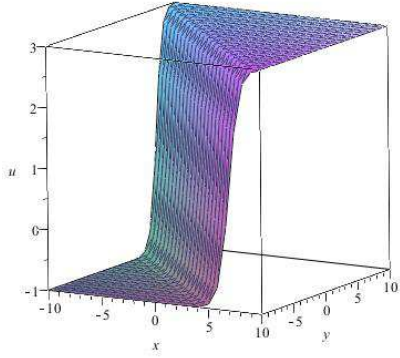}\label{sol15u}}
  \hfill
  \subfloat[]{\includegraphics[scale=0.30]{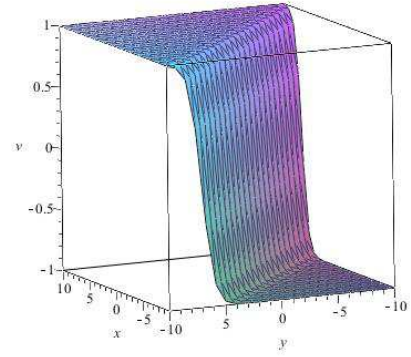}\label{sol15v}}
  \caption{Solution profile of \eqref{gov} for a solution \eqref{sol15}: (a) 3d profile of $u_{15}$ at $t=1$ (b) 3d profile of $v_{15}$ at $t=1$}
  \label{fsol15}
  \end{figure}

which yields another exact solution for the given system \eqref{gov} as follows:
\begin{eqnarray}\label{sol15}
\nonumber&&u_{15}=\frac{k_3\cosh(C_{35}(k_1x+k_2y+k_3t+k_4)+C_{36})+2ck_1^2C_{35}\sinh(C_{35}(k_1x+k_2y+k_3t+k_4)+C_{36})}
{gk_1\cosh(C_{35}(k_1x+k_2y+k_3t+k_4)+C_{36})},\\
&&v_{15}=\tanh(C_{35}(k_1x+k_2y+k_3t+k_4)+C_{36})
\end{eqnarray}
where $C_{35}$ and $C_{36}$ are arbitrary constants provided that $a=\frac{c^2k_2C_{35}}{b-ck_2C_{35}}$ and $g=\frac{2c^2k_2C_{35}}{b-ck_2C_{35}}$.

Now we demonstrate the physical significance of the solution profile for $u_{15}$ and $v_{15}$ given by \eqref{sol15} in the Figure \ref{fsol15} with respect to $x$ and $y$ at fixed $t=1$ by considering the parameters $b=3,c=1,k_1=1,k_2=1,k_3=1,k_4=1,C_{35}=1,C_{36}=1,a=\frac{c^2k_2C_{35}}{b-ck_2C_{35}}=\frac{1}{2}$ and $g=\frac{2c^2k_2C_{35}}{b-ck_2C_{35}}=1$. Here we noticed that the solution profile for $u_{15}$, illustrated in the Figure \ref{sol15u}, represents a kink type soliton whilst the solution profile for $v_{15}$ displayed in the Figure \ref{sol15v} represents an anti-kink type soliton profile.


Now, by imposing the condition $a=\frac{c}{2}$ and $g=c$, we have another exact solution of the above reduced system of ODEs \eqref{red6} as
\begin{eqnarray*}
&&U(h)=\frac{k_3C_{38}-2ck_1^2C_{37}+2k_3C_{37}h}{ck_1(C_{38}+2C_{37}h)},\\
&&V(h)=C_{37}h^2+C_{38}h+C_{39}
\end{eqnarray*}

\begin{figure}[h!]
  \centering
  \subfloat[]{\includegraphics[scale=0.28]{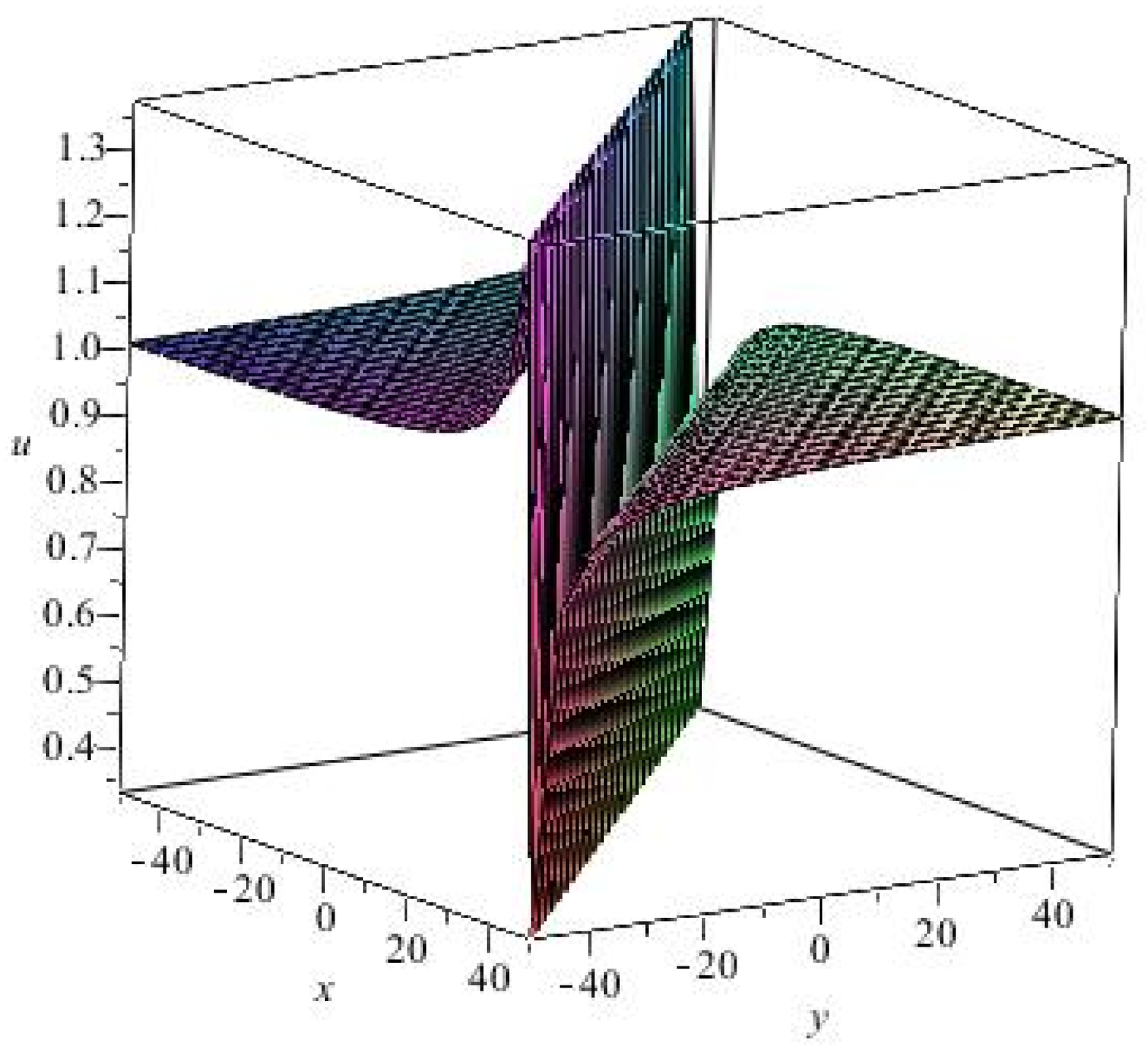}\label{sol16u1}}
  \subfloat[]{\includegraphics[scale=0.28]{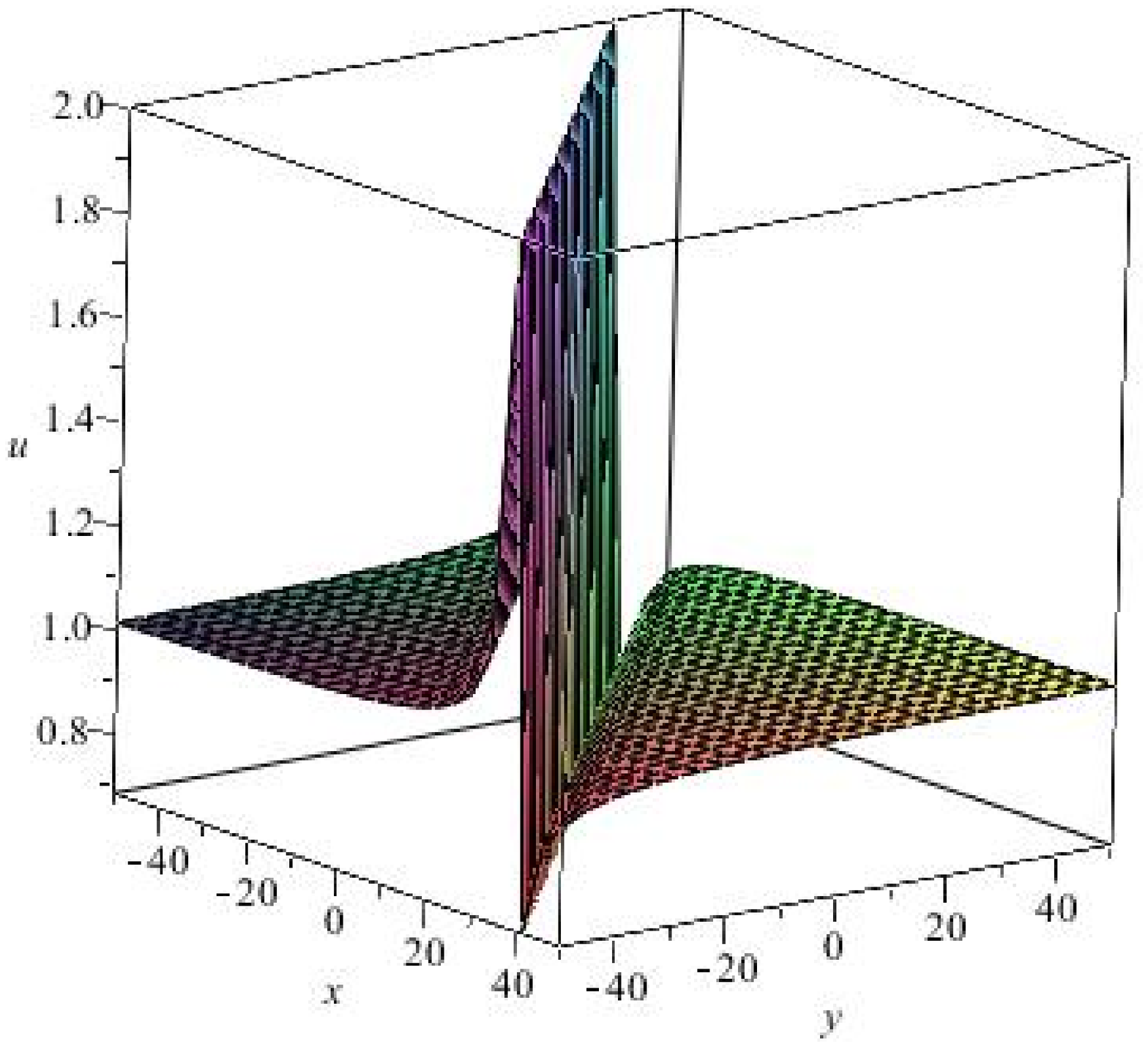}\label{sol16u2}}
  \subfloat[]{\includegraphics[scale=0.28]{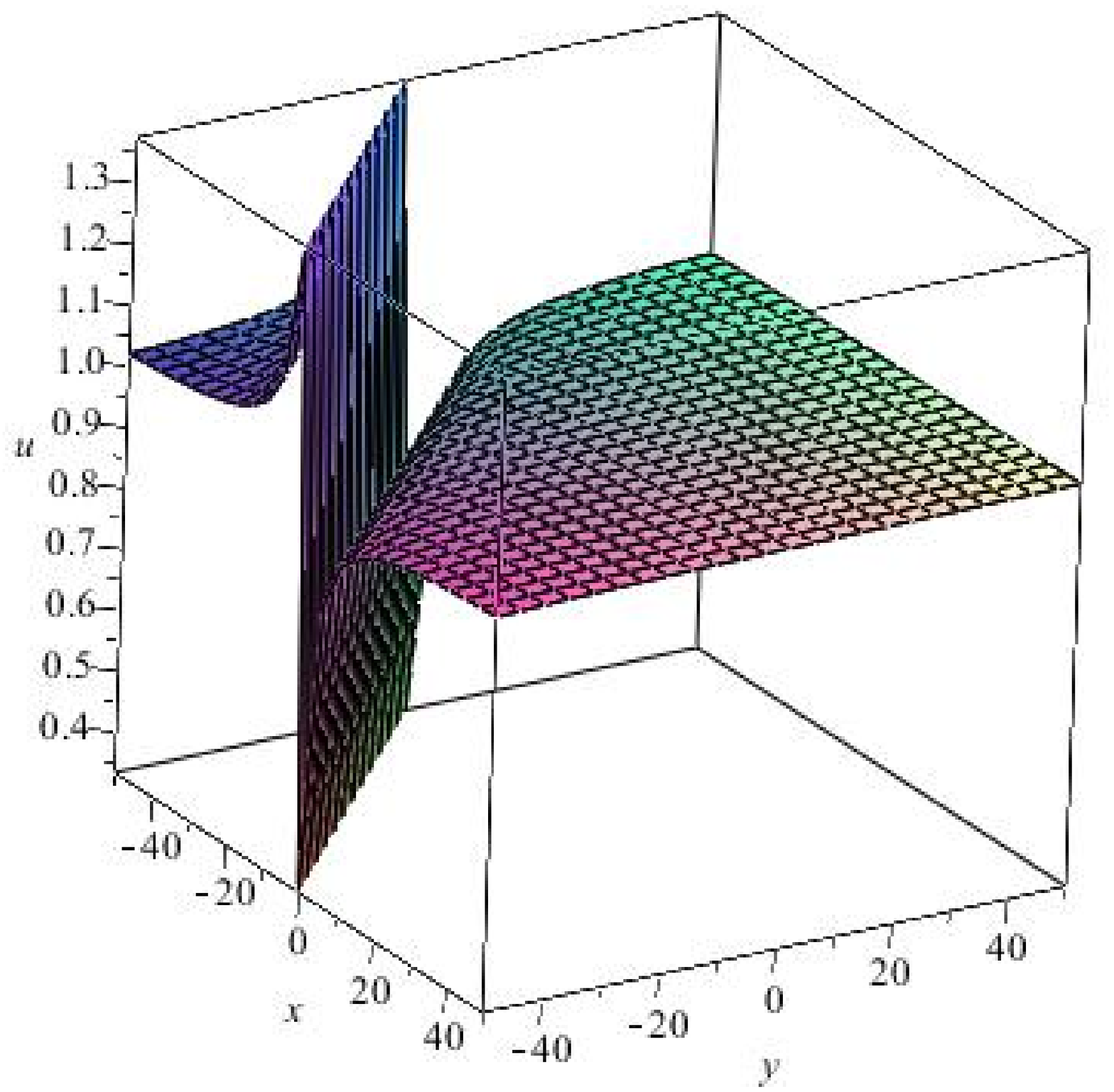}\label{sol16u3}}

  \subfloat[]{\includegraphics[scale=0.28]{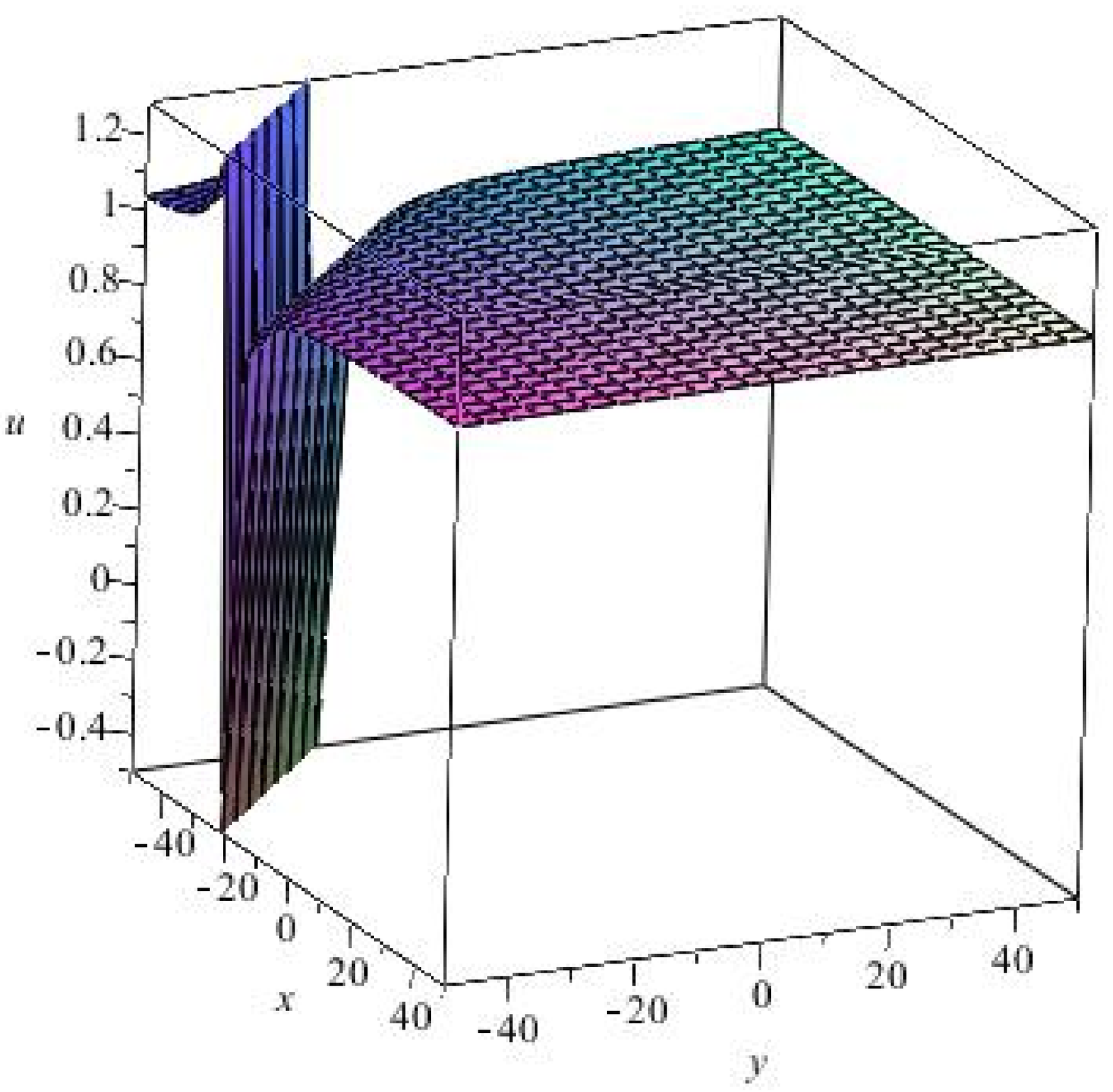}\label{sol16u4}}
  \subfloat[]{\includegraphics[scale=0.28]{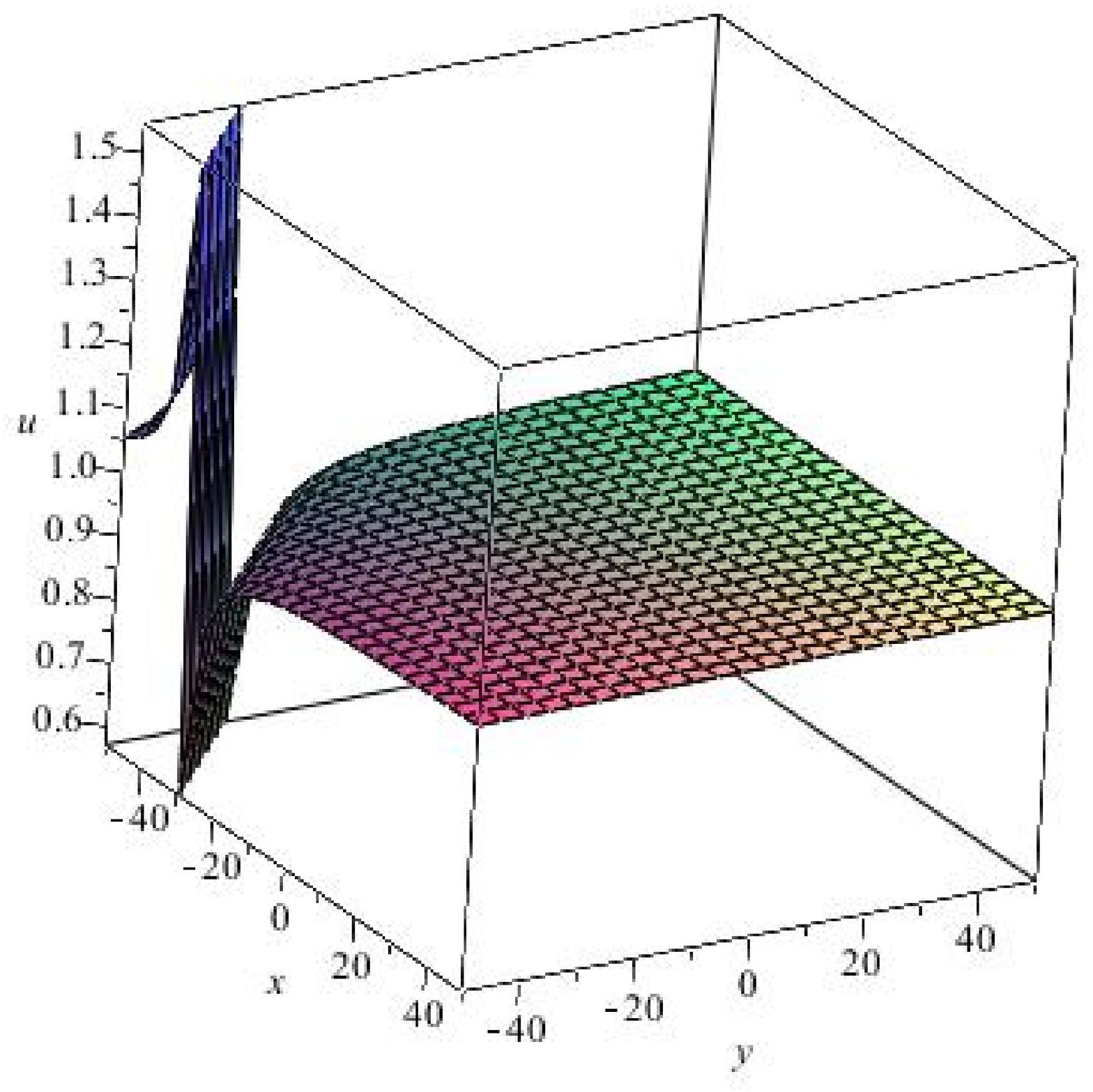}\label{sol16u5}}
  \subfloat[]{\includegraphics[scale=0.28]{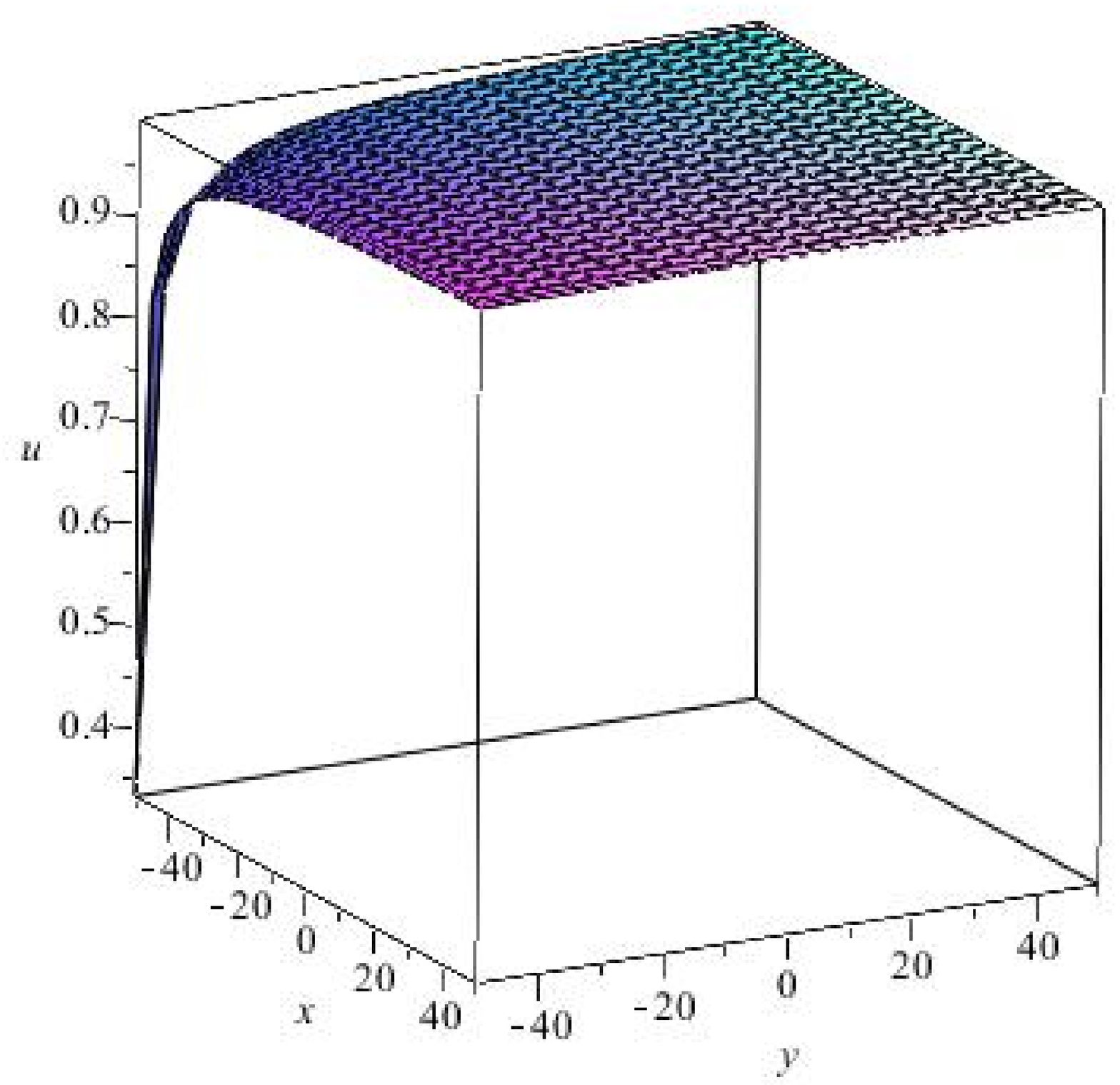}\label{sol16u6}}

  \subfloat[]{\includegraphics[scale=0.28]{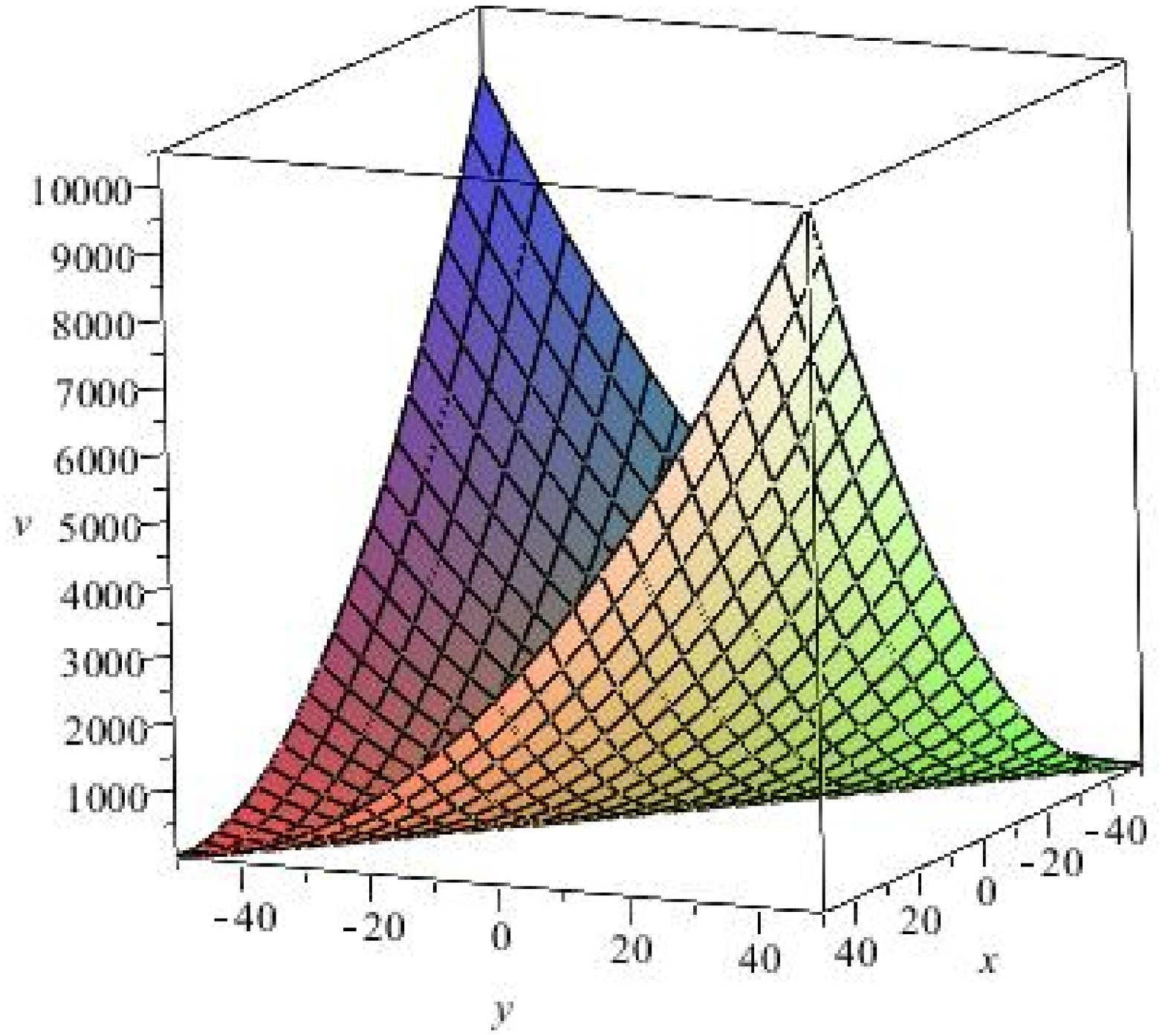}\label{sol16v}}
  \caption{Solution profile for a solution \eqref{sol16}: (a) $u_{16}$ at $t=0$, (b) $u_{16}$ at $t=10$, (c) $u_{16}$ at $t=50$, (d) $u_{16}$ at $t=70$, (e) $u_{16}$ at $t=80$, (f) $u_{16}$ at $t=100$, (g) $v_{16}$ at $t=1$}
  \label{fsol16}
  \end{figure}

which results another exact solution of the given system \eqref{gov} of the form as follows:
\begin{eqnarray}\label{sol16}
&&u_{16}=\frac{k_3C_{38}-2ck_1^2C_{37}+2k_3C_{37}(k_1x+k_2y+k_3t+k_4)}{ck_1(C_{38}+2C_{37}(k_1x+k_2y+k_3t+k_4))},\\
\nonumber&&v_{16}=C_{37}(k_1x+k_2y+k_3t+k_4)^2+C_{38}(k_1x+k_2y+k_3t+k_4)+C_{39}
\end{eqnarray}
where $C_{37},C_{38}$ and $C_{39}$ are arbitrary constants.

We investigate the physical behavior of the solution profile of \eqref{sol16} for $u_{16}$ and $v_{16}$ in the Figure \ref{fsol16} by choosing $c=1,k_1=1,k_2=1,k_3=1,k_4=1,C_{37}=1,C_{38}=1$ and $C_{39}=1$. We noticed from the Figure \ref{fsol16} that $u_{16}$ represents a singular kink traveling wave solution profile for different values of $t$. It is very interesting to observe that as time evolves, the singular kink waves (see, Figure \ref{sol16u1}-\ref{sol16u6}) travel towards the negative $x$-direction and after  certain time the kink property disappears. On the other hand, we demonstrate the physical behavior of $v_{16}$ at $t=1$ in the Figure \ref{sol16v} and observed  bowl shaped profile.

\section{Conservation laws}\label{s6}
Conservation laws deal with essential physical properties of the process modeled by a given PDE system and have also wide applications in existence, uniqueness and stability analysis for the development of numerical methods. Moreover, one can construct a nonlocally related PDE systems of the original PDE system by
introducing some potential (nonlocal) variables through conservation laws and thus possibility of finding nonlocal symmetries and hence new exact solutions. A conservation law of the given PDE system \eqref{gov} can be represented in divergence form as $D_t\phi^t+D_x\phi^x+D_y\phi^y=0$ which holds true for the solution manifold of the same system. Recently, Anco and Bluman \cite{anco1997direct,anco2002direct,anco2002direct1} presented a systematic procedure to construct conservation law multipliers using direct multiplier method in terms of Euler operator which annihilates divergence expression. The advantage of this method over that of Noether's method is that this method does not require that the given PDE system to admit variational symmetry. In \cite{wang20202+} the authors applied  direct multipliers method to obtain conservation law multipliers by considering $a = 1, b = 2, c = 1$ and $g = 2$, consequently constructed the corresponding conserved vectors. Authors claimed that $$Q=f(y)v+g_1(t)x^2+g_2(t)x+F(y)+g_3(t)$$ is the multiplier of \eqref{govv} where $f(y),g_1(t),g_2(t),g_3(t)$ and $F(y)$ are arbitrary functions. This claim lacks proper sense of understanding as there should be a set of multipliers where each set consists of two quantities instead of just one quantity. As a result, the associated conserved vectors are also incorrect. \\
So, in this section, we apply the direct multiplier technique to the original system \eqref{gov} and construct conservation laws systematically. The Euler operator with respect to dependent variables $u^j,j=1,2$ $(u^1=u,u^2=v)$, is the operator defined by
$$E_{u^j}=\frac{\partial}{\partial u^j}-D_i\frac{\partial}{\partial u^j_i}+...+(-1)^lD_{i_1}...D_{i_l}\frac{\partial}{\partial u^j_{i_1...i_l}}+...$$
 So, using the multipliers method we obtain the following sets of conservation law multipliers
\begin{eqnarray*}
&\Lambda_1^1=\frac{1}{2}\alpha(t)x^2,&\Lambda_2^1=0,\\
&\Lambda_1^2=\beta(t)x,&\Lambda_2^2=0,\\
&\Lambda_1^3=\gamma(t),&\Lambda_2^3=0,\\
&\Lambda_1^4=\mu(y),&\Lambda_2^4=0
\end{eqnarray*}
where $\alpha(t),\beta(t),\gamma(t)$ and $\mu(y)$ are arbitrary functions. Using these multipliers we obtain the conservation laws of the form $D_t\phi_j^t+D_x\phi_j^x+D_y\phi_j^y=0,~~j=1,2,3,4$ whose conserved vectors are given by
\begin{eqnarray}
&&\phi_1^t=\frac{1}{2}x^2u_y\alpha(t),\phi_1^x=\left(-auu_y-axu_y+\frac{1}{2}ax^2u_{xy}-bv+bxv_x-\frac{1}{2}bx^2v_{xx}\right)\alpha(t),\\
\nonumber&&\phi_1^y=-\frac{1}{2}u(-2axu\alpha(t)-2a\alpha(t)+x^2\alpha'(t));\\ \nonumber&&\\
&&\phi_2^t=xu_y\beta(t),\phi_2^x=(-2axuu_y-au_y+axu_{xy}+bv_x-bv_{xx})\beta(t),\\
\nonumber&&\phi_2^y=-u(-au\beta(t)+x\beta'(t));\\ \nonumber&&\\
&&\phi_3^t=u_y\gamma(t),\phi_3^x=(-2auu_y+au_{xy}-bv_{xx})\gamma(t),\phi_3^y=-u\gamma'(t)~~\text{and}\\ \nonumber&&\\
&&\phi_4^t=u_y\mu(y),\phi_4^x=(-2auu_y+au_{xy}-bv_{xx})\mu(y),\phi_4^y=0.
\end{eqnarray}
Conservation laws are useful for various applications including construction of nonlocally related PDE systems. Also, one can perform nonlocal symmetry analysis and further construct nonlocal conservation laws those are very challenging and recent topic of research. We show the existence of nonlocal conservation laws of the given system \eqref{gov} in the succeeding section.

\section{Applications}\label{s7}
Conservation laws are very useful in constructing nonlocally related PDE systems, developing mathematical theory of nonlocal conservation laws and nonlocal symmetry analysis and thus new exact solutions. It would be very interesting and challenging to perform nonlocal symmetry analysis of (2+1)-dimensional nonlinear system of PDEs. This kind of problems are stated as open problems in \cite{cheviakov2010multidimensional,cheviakov2010multidimensional2}. For PDE systems with $n>3$ independent variables, the situation for obtaining and using nonlocally related PDE systems is considerably more complex than in the case of $n=2$. In particular, every divergence-type
conservation law gives rise to several potential variables, which are only defined to within arbitrary functions of the independent variables. The corresponding potential system is thus under-determined, and is said to have gauge freedom. Additional equations involving potential variables, called gauge constraints, are needed to make such potential systems determined.

For example,
consider a divergence-type conservation law in three-dimensional space
$$\text{div} \phi=\phi^1_x+\phi^2_y+\phi^3_z=0$$

with flux vector $\phi=(\phi^1(x,y,z),\phi^2(x,y,z),\phi^3(x,y,z))$ and independent variables $x,y,z$. It immediately follows that there exists a vector potential $\psi=(\psi^1(x,y,z),\psi^2(x,y,z),\psi^3(x,y,z))$, such that $\phi=\text{curl} \psi$. Consequently, the potential system in this case becomes
\begin{eqnarray*}
&&\psi^3_y-\psi^2_z=\phi^1,\\
&&\psi^1_z-\psi^3_x=\phi^2,\\
&&\psi^2_x-\psi^1_y=\phi^3.
\end{eqnarray*}
However, unlike in the two-dimensional situation, the   potential
system  is under-determined. An additional equation involving the potential variables is required in order
to complete the potential
system to eliminate its gauge
freedom. For example, one can have the gauges:\\\\
$\bullet$ divergence (Coulomb) gauge: div $\psi=\psi^1_x+\psi^2_y+\psi^3_z=0,$\\
$\bullet$ spatial gauge: $\psi^k=0,$ $k=1$ or 2 or 3,\\
$\bullet$ Poincar$\acute{e}$ gauge: $x\psi^1+y\psi^2+z\psi^3=0$,\\
provided that all solutions of the potential system can be obtained from
the solution of the corresponding gauge-constrained (determined) potential
system. If one of the coordinates in a given PDE system is time $t$, special gauges
are frequently used, such as\\\\
$\bullet$ Lorentz gauge (in (2+1)-dimensional): $\psi^1_t-\psi^2_x-\psi^3_y=0$,\\
$\bullet$ Cronstrom gauge (in (2+1)-dimensional): $t\psi^1-x\psi^2-y\psi^3=0.$\par
Here we write down the potential systems associated to the given system \eqref{gov} by making use of the conserved vectors given in preceding section.
\begin{eqnarray}\label{pot1}
\text{Potential system I:}&&\psi^3_x-\psi^2_y=\frac{1}{2}x^2u_y\alpha(t),\\
\nonumber&&\psi^1_y-\psi^3_t=\left(-auu_y-axu_y+\frac{1}{2}ax^2u_{xy}-bv+bxv_x-\frac{1}{2}bx^2v_{xx}\right)\alpha(t),\\
\nonumber&&\psi^2_t-\psi^1_x=-\frac{1}{2}u(-2axu\alpha(t)-2a\alpha(t)+x^2\alpha'(t));
\end{eqnarray}

\begin{eqnarray}\label{pot2}
\text{Potential system II:}&&\psi^3_x-\psi^2_y=xu_y\beta(t),\\
\nonumber&&\psi^1_y-\psi^3_t=(-2axuu_y-au_y+axu_{xy}+bv_x-bv_{xx})\beta(t),\\
\nonumber&&\psi^2_t-\psi^1_x=-u(-au\beta(t)+x\beta'(t));
\end{eqnarray}

\begin{eqnarray}\label{pot3}
\text{Potential system III:}&&\psi^3_x-\psi^2_y=u_y\gamma(t),\\
\nonumber&&\psi^1_y-\psi^3_t=(-2auu_y+au_{xy}-bv_{xx})\gamma(t),\\
\nonumber&&\psi^2_t-\psi^1_x=-u\gamma'(t);
\end{eqnarray}
and

\begin{eqnarray}\label{pot4}
\text{Potential system IV:}&&\psi^3_x-\psi^2_y=u_y\mu(y),\\
\nonumber&&\psi^1_y-\psi^3_t=(-2auu_y+au_{xy}-bv_{xx})\mu(y),\\
\nonumber&&\psi^2_t-\psi^1_x=0.
\end{eqnarray}

For the potential system I, let us consider the spatial gauge $\psi^3=0$ given in \eqref{pot1} and by choosing $\alpha(t)=2$, the corresponding conservation law multipliers are
\begin{eqnarray}\label{mult1}
\nonumber&&\delta^1=\left(\frac{\psi^2}{x^2}+u\right)k_1(t),\\
&&\delta^2=\frac{H_1(t)x+H_2(t)}{x^3},\\
\nonumber&&\delta^3=0,~~\delta^4=0
\end{eqnarray}

and

\begin{eqnarray}\label{mult2}
\nonumber&&\delta^5=\left(\frac{\psi^2}{x^2}+u\right)k_2(t),\\
&&\delta^6=\frac{H_3(t)x+H_4(t)H_5(y)+H_6(t)}{x^3},\\
\nonumber&&\delta^7=-\frac{1}{2}\frac{H_4(t)H_5'(y)}{x^2},~~\delta^8=0
\end{eqnarray}
where $k_1(t),k_2(t),H_1(t),H_2(t),H_3(t),H_4(t)$ and $H_5(y)$ are arbitrary functions.
Here we observe that the multiplier components $\delta^1$ and $\delta^5$ have an essential dependence on the potential variable $\psi^2.$ Hence they yield nonlocal conservation laws of the given system \eqref{gov}. Similarly, it is easy to show that potential system II and potential system IV also yield nonlocal conservation laws for \eqref{gov} by considering $\beta(t)=1$ and $\mu(y)=1.$ \par
In future study, it will be very interesting to analyze nonlocal symmetries of the given system \eqref{gov} arising from potential systems (I-IV) as well as from inverse potential systems and to obtain some new exact solutions.

\section{Conclusions}\label{s8}
The (2+1)-dimensional BLP system is studied in the context of classical Lie symmetry analysis. It is observed that the system admits infinite dimensional Lie algebra. We performed the classification of optimal subalgebras and using each subalgebra we obtained several new exact solutions (\eqref{sol1}, \eqref{sol2}, \eqref{sol3}, \eqref{sol4}, \eqref{sol5}, \eqref{sol6}, \eqref{sol7}, \eqref{sol8}, \eqref{sol9}, \eqref{sol10}, \eqref{sol11}, \eqref{sol12}, \eqref{sol13}, \eqref{sol14}, \eqref{sol15} and \eqref{sol16}) of the BLP system. In addition to that, we noticed that the only solution \eqref{recover} presented  in \cite{wang20202+} was recovered as a particular case of the obtained solution \eqref{sol12}. The computed solutions are reported first time in the literature. Physical behavior of some of the solutions are exhibited geometrically with the help of numerical simulations which consists of traveling waves, lump type solitons, kink and anti-kink type solitons, breather solitons, singular kink type solitons and etc. We constructed some conservation laws of the given system by using the direct multipliers method those may be used to further study on the nonlocal symmetry analysis of BLP system. Finally, as an application, we study the nonlocal conservation laws of the given system by using direct multiplier technique to the corresponding potential systems and appending spatial gauge constraints on them.\vskip6pt

\section*{Acknowledgement}
First author is highly thankful to Ministry of Human Resource Development, Government of India, for the institute fellowship (grant no.IIT/ACAD/PGS$\&$R/F.II/2/16MA90J04).

\section*{Conflict of interest}
The authors declare that they have no conflict of interest.

\section*{References}
\bibliography{ref}

\begin{thebibliography}{10}
\expandafter\ifx\csname url\endcsname\relax
  \def\url#1{\texttt{#1}}\fi
\expandafter\ifx\csname urlprefix\endcsname\relax\def\urlprefix{URL }\fi
\expandafter\ifx\csname href\endcsname\relax
  \def\href#1#2{#2} \def\path#1{#1}\fi

\bibitem{bluman2010applications}
G.~W. Bluman, A.~F. Cheviakov, S.~C. Anco, Applications of symmetry methods to
  partial differential equations, Vol. 168, Springer, 2010.

\bibitem{olver1987group}
P.~J. Olver, P.~Rosenau, Group-invariant solutions of differential equations,
  SIAM Journal on Applied Mathematics 47~(2) (1987) 263--278.

\bibitem{sahoo2020optimal}
S.~M. Sahoo, T.~Raja~Sekhar, G.~P. Raja~Sekhar, Optimal classification, exact
  solutions, and wave interactions of euler system with large friction,
  Mathematical Methods in the Applied Sciences 43~(9) (2020) 5744--5757.

\bibitem{satapathy2018optimal}
P.~Satapathy, T.~Raja~Sekhar, Optimal system, invariant solutions and evolution
  of weak discontinuity for isentropic drift flux model, Applied Mathematics
  and Computation 334 (2018) 107--116.

\bibitem{yacsar2011symmetries}
E.~Ya{\c{s}}ar, T.~{\"O}zer, On symmetries, conservation laws and invariant
  solutions of the foam-drainage equation, International Journal of Non-Linear
  Mechanics 46~(2) (2011) 357--362.

\bibitem{saha2020lie2}
S.~Saha~Ray, Vinita, Lie symmetry analysis, symmetry reductions with exact
  solutions, and conservation laws of (2+ 1)-dimensional bogoyavlenskii-schieff
  equation of higher order in plasma physics, Mathematical Methods in the
  Applied Sciences 43~(9) (2020) 5850--5859.

\bibitem{ovsiannikov2014group}
L.~V. Ovsiannikov, Group analysis of differential equations, Academic press,
  2014.

\bibitem{olver2000applications}
P.~J. Olver, Applications of Lie groups to differential equations, Vol. 107,
  Springer Science \& Business Media, 2000.

\bibitem{sekhar2016group}
T.~Raja~Sekhar, P.~Satapathy, Group classification for isothermal drift flux
  model of two phase flows, Computers \& Mathematics with Applications 72~(5)
  (2016) 1436--1443.

\bibitem{cherniha2021complete}
R.~Cherniha, M.~Serov, Y.~Prystavka, A complete lie symmetry classification of
  a class of (1+ 2)-dimensional reaction-diffusion-convection equations,
  Communications in Nonlinear Science and Numerical Simulation 92 (2021)
  105466.

\bibitem{sil2020nonclassical}
S.~Sil, T.~Raja~Sekhar, Nonclassical symmetry analysis, conservation laws of
  one-dimensional macroscopic production model and evolution of nonlinear
  waves, Journal of Mathematical Analysis and Applications (2020) 124847.

\bibitem{vaneeva2020generalization}
O.~O. Vaneeva, A.~Bihlo, R.~O. Popovych, Generalization of the algebraic method
  of group classification with application to nonlinear wave and elliptic
  equations, Communications in Nonlinear Science and Numerical Simulation 91
  (2020) 105419.

\bibitem{benoudina2021lie}
N.~Benoudina, Y.~Zhang, C.~M. Khalique, Lie symmetry analysis, optimal system,
  new solitary wave solutions and conservation laws of the pavlov equation,
  Communications in Nonlinear Science and Numerical Simulation 94 (2021)
  105560.

\bibitem{opanasenko2017group}
S.~Opanasenko, A.~Bihlo, R.~O. Popovych, Group analysis of general
  burgers--korteweg--de vries equations, Journal of Mathematical Physics 58~(8)
  (2017) 081511.

\bibitem{zhang2021symmetry}
Z.-Y. Zhang, G.-F. Li, Symmetry properties of conservation laws for nonlinear
  fokker-planck equation describing cell population growth, Communications in
  Nonlinear Science and Numerical Simulation 93 (2021) 105506.

\bibitem{liu2020existence}
M.~Liu, H.~Dong, On the existence of solution, lie symmetry analysis and
  conservation law of magnetohydrodynamic equations, Communications in
  Nonlinear Science and Numerical Simulation 87 (2020) 105277.

\bibitem{anco1997direct}
S.~C. Anco, G.~W. Bluman, Direct construction of conservation laws from field
  equations, Physical Review Letters 78~(15) (1997) 2869.

\bibitem{anco2002direct}
S.~C. Anco, G.~W. Bluman, Direct construction method for conservation laws of
  partial differential equations part i: Examples of conservation law
  classifications, European Journal of Applied Mathematics 13~(5) (2002)
  545--566.

\bibitem{anco2002direct1}
S.~C. Anco, G.~W. Bluman, Direct construction method for conservation laws of
  partial differential equations part ii: General treatment, European Journal
  of Applied Mathematics 13~(5) (2002) 567--585.

\bibitem{sil2020nonlocal}
S.~Sil, T.~Raja~Sekhar, D.~Zeidan, Nonlocal conservation laws, nonlocal
  symmetries and exact solutions of an integrable soliton equation, Chaos,
  Solitons \& Fractals 139 (2020) 110010.

\bibitem{sil2020nonlocally}
S.~Sil, T.~Raja~Sekhar, Nonlocally related systems, nonlocal symmetry
  reductions and exact solutions for one-dimensional macroscopic production
  model, The European Physical Journal Plus 135~(6) (2020) 1--23.

\bibitem{BLP}
M.~Boiti, J.~J.-P. Leon, F.~Pempinelli, Integrable two-dimensional
  generalisation of the sine- and sinh-gordon equations, Inverse Problems 3~(1)
  (1987) 37--49.

\bibitem{mu2013localized}
G.~Mu, Z.~Dai, Z.~Zhao, Localized structures for (2+ 1)-dimensional
  boiti--leon--pempinelli equation, Pramana 81~(3) (2013) 367--376.

\bibitem{jiang2010solitons}
Y.~Jiang, B.~Tian, W.-J. Liu, M.~Li, P.~Wang, K.~Sun, Solitons, b{\"a}cklund
  transformation, and lax pair for the (2+ 1)-dimensional
  boiti--leon--pempinelli equation for the water waves, Journal of mathematical
  physics 51~(9) (2010) 093519.

\bibitem{zhao2017lie}
Z.~Zhao, B.~Han, Lie symmetry analysis, b{\"a}cklund transformations, and exact
  solutions of a (2+ 1)-dimensional boiti-leon-pempinelli system, Journal of
  Mathematical Physics 58~(10) (2017) 101514.

\bibitem{lu2004soliton}
Z.~L{\"u}, H.~Zhang, Soliton like and multi-soliton like solutions for the
  boiti--leon--pempinelli equation, Chaos, Solitons \& Fractals 19~(3) (2004)
  527--531.

\bibitem{ma2003diversity}
W.-X. Ma, Diversity of exact solutions to a restricted boiti--leon--pempinelli
  dispersive long-wave system, Physics Letters A 319~(3-4) (2003) 325--333.

\bibitem{huang2004exact}
D.-J. Huang, H.-Q. Zhang, Exact travelling wave solutions for the
  boiti--leon--pempinelli equation, Chaos, Solitons \& Fractals 22~(1) (2004)
  243--247.

\bibitem{wazwaz2010}
A.-M. Wazwaz, M.~S. Mehanna, A variety of exact travelling wave solutions for
  the (2 + 1)-dimensional boiti-leon-pempinelli equation, Communications in
  Theoratical Physics 217.

\bibitem{wang2014}
Y.-H. Wang, H.~Wang, Symmetry analysis and cte solvability for the
  (2+1)-dimensional boiti-leon-pempinelli equation, Physica Scripta 89.

\bibitem{yang2011}
Y.~Zheng, M.~Song-Hua, F.~Jian-Ping, Soliton excitations and chaotic patterns
  for the (2+1)-dimensional boiti- leon-pempinelli system, Chinese Physics B
  20~(6) (2011) 060506.

\bibitem{kumar2014}
M.~Kumar, R.~Kumar, On new similarity solutions of the boiti-leon-pempinelli
  system, Communications in Theoratical Physics 61~(1) (2014) 121--126.

\bibitem{wang20202+}
G.~Wang, J.~Vega-Guzman, A.~Biswas, A.~K. Alzahrani, A.~H. Kara, (2+
  1)-dimensional boiti--leon--pempinelli equation--domain walls, invariance
  properties and conservation laws, Physics Letters A 384~(10) (2020) 126255.

\bibitem{Malek2018Lie}
M.~B. Abd-el Malek, A.~M. Amin, Lie group method for solving viscous barotropic
  vorticity equation in ocean climate models, Computers and Mathematics with
  Applications 75.

\bibitem{cheviakov2010multidimensional}
A.~F. Cheviakov, G.~W. Bluman, Multidimensional partial differential equation
  systems: Generating new systems via conservation laws, potentials, gauges,
  subsystems, Journal of mathematical physics 51~(10) (2010) 103521.

\bibitem{cheviakov2010multidimensional2}
A.~F. Cheviakov, G.~W. Bluman, Multidimensional partial differential equation
  systems: Nonlocal symmetries, nonlocal conservation laws, exact solutions,
  Journal of mathematical physics 51~(10) (2010) 103522.

\end{thebibliography}
\bibliographystyle{unsrt}


\end{document}